\documentclass[12pt]{article}
\usepackage{amsmath,amsthm,amsfonts,amssymb}

\title{\rt s in topology, geometry, and group theory}
\author{Mladen Bestvina\thanks{Support by the National Science
Foundation is gratefully acknowledged.}}
\date{\today\\ preliminary version}

\newtheorem{thm}{Theorem}[section]
\newtheorem{lemma}[thm]{Lemma}

\newtheorem{prop}[thm]{Proposition}

\theoremstyle{remark}
\newtheorem{example}[thm]{Example}
\newtheorem{exercise}[thm]{Exercise}
\newtheorem{definition}[thm]{Definition}

\def\C{{\mathbb C}}
\def\R{{\mathbb R}}
\def\Z{{\mathbb Z}}
\def\H{{\mathbb H}}
\def\rt{$\R$-tree}
\def\dh{$\delta$-hyperbolic}

\begin{document}

\maketitle
\tableofcontents

\section{Introduction}

This paper is intended as a survey of the theory and applications of
real trees from a topologist's point of view. The idea of an
all-inclusive historical account was quickly abandoned at the start of
this undertaking, but I hope to describe the main ideas in the
subject with emphasis on applications outside the theory of
$\R$-trees. The ``Rips machine'', i.e. the classification of measured
laminations on 2-complexes, is the key ingredient. Roughly speaking,
the Rips machine is an algorithm that takes as input a finite
2-complex equipped with a transversely measured lamination (more
precisely, a band complex), and puts it in a ``normal form''. This
normal form is surprisingly simple -- the lamination is the disjoint
union of finitely many sub-laminations each of which belongs to one of
four types:

$\bullet$ {\it simplicial}: all leaves are compact and the lamination
is a bundle over a leaf with compact 0-dimensional fiber,

$\bullet$ {\it surface}: geodesic lamination on a compact hyperbolic
surface (or a cone-type orbifold),

$\bullet$ {\it toral}: start with a standard lamination of the $n$-torus
by irrational planes of codimension 1 and restrict to the 2-skeleton; more
generally, replace the torus by a cone-type orbifold covered by a torus
(with the deck group leaving the lamination invariant),

$\bullet$ {\it thin}: this type is most interesting of all. It was
discovered and studied by G. Levitt \cite{gl:thin}. See section 5.3 for
the definition and basic properties.

\vskip .2cm Measured laminations on 2-complexes arise in the study of
$\R$-trees via a process called {\it resolution}. In the simplicial
case, this idea goes back to J. Stallings and was used with great
success by M. Dunwoody.  If $G$ is a finitely presented group that
acts by isometries on an $\R$-tree, one wants to deduce the structure
of $G$, given the knowledge of vertex and arc stabilizers. Bass-Serre
theory \cite{se:trees} solves this beautifully in the case of
simplicial trees. For an exposition of Bass-Serre theory from a
topological point of view, see \cite{sw:topmethods}.

I hope to convince the reader that the development of the theory of
$\R$-trees is not an idle exercise in generalizations -- indeed, in
addition to the intrinsic beauty of the theory, $\R$-trees appear in
``real life'' as a brief look at the final section of this survey
reveals. The reason for this is the construction presented in section
3, which takes a sequence of isometric actions of $G$ on ``negatively
curved spaces'' and produces an isometric action of $G$ on an
$\R$-tree in the (Gromov-Hausdorff) limit.

The central part of the paper (sections 4-6) is devoted to a study of
the Rips machine and the structure theory of groups that act
isometrically on $\R$-trees. The approach follows closely
\cite{bf:stable}, and the reader is referred to that paper for more
details.  Gaboriau, Levitt, and Paulin have developed a different (but
equivalent) point of view in a series of papers (see references, and
in particular the survey \cite{paulin:bourbaki} which puts everything
together).  For the historical developments and the state of the
theory preceding Rips' breakthrough, see the surveys
\cite{ps:dendrology1} and \cite{ps:dendrology2}.

This paper is an expanded version of a talk presented at the AMS
meeting \# 906 in Greensboro, NC in October 1995. I would like to
thank Mark Feighn, Gilbert Levitt, and Zlil Sela for useful comments.
I will always be grateful to Mark Feighn for our long-term collaboration
and for all the fun we had while learning and contributing to the
mathematics described here.

\section{Definition and first examples of \rt s}

\begin{definition}
Let $(X,d)$ be a metric space and let $x,y\in X$. An {\it arc from $x$ to $y$}
is the image of a topological embedding $\alpha:[a,b]\to X$ of a closed
interval $[a,b]$ (and we allow the possibility $a=b$) such that
$\alpha(a)=x$ and $\alpha(b)=y$. A {\it geodesic segment} from $x$
to $y$ is the image of an isometric
embedding $\alpha:[a,b]\to X$ with $\alpha(a)=x$ and $\alpha(b)=y$.
\end{definition}

\begin{definition}
We say that $(X,d)$ is an {\it \rt} if for any $x,y\in X$ there is a unique
arc from $x$ to $y$ and this arc is a geodesic segment.
\end{definition}

\begin{example}
Let $X$ be a connected 1-dimensional simplicial complex that contains no
circles. For every edge $e$ of $X$ choose an embedding $e\to\R$. If $x,y\in X$,
there is a unique arc $A$ from $x$ to $y$. This arc can be subdivided into
subarcs $A_1,A_2,\cdots,A_n$ each of which is contained in an edge of $X$.
Define the length of $A_i$ as the length of its image in $\R$ under the
chosen embedding, and define $d(x,y)$ as the sum of the lengths of
the $A_i$'s. The metric space $(X,d)$ is an \rt. We say that an \rt\ is
{\it simplicial} if it arises in this fashion.
\end{example}

\begin{example} (SNCF metric)
Take $X=\R^2$ and let $e$ denote the Euclidean distance on $X$.
Define a new distance $d$ as follows. We imagine that there is a train line
operating along each ray from the origin (=Paris). If two points $x,y\in X$
lie on the same ray, then $d(x,y)=e(x,y)$. In all other cases the train ride
from $x$ to $y$ goes through the origin, so $d(x,y)=e(0,x)+e(0,y)$.
The metric space $(X,d)$ is a (simplicial) \rt.
\end{example}

\begin{example}
A slight modification of the previous example yields a non-simplicial \rt.
Take $X=\R^2$ and imagine trains operating on all vertical lines as well as
along the $x$-axis. Thus $d(x,y)=e(x,y)$ when $x,y$ are on the same
vertical line, and $d(x,y)=|x_2|+|y_2|+|x_1-y_1|$ otherwise, where we
set $x=(x_1,x_2)$ and $y=(y_1,y_2)$. 
\end{example}

$\R$-trees that arise in applications tend to be separable (as metric spaces),
and in fact they are the union of countably many lines. Example 2.5 can be
easily modified to yield an example of a separable non-simplicial $\R$-tree
(restrict to the subset of $X$ consisting of the $x$-axis and the points
with rational $x$-coordinate).

\subsection{Isometries of \rt s}

I will now recall basic facts about isometric actions on
$\R$-trees. Proofs are a straightforward generalization from the case
of simplicial trees that can be found in
\cite{se:trees}. Alternatively, the reader is referred to
\cite{ms:valuations1}, \cite{cm:trees}, or \cite{ab:lengthfunctions}.

Let $\phi:T\to T$ be an isometry of an $\R$-tree $T$. The {\it translation
length} of $\phi$ is the number
$$\ell(\phi)=\inf \{ d(x,\phi(x))|x\in T\}$$ where $d$ denotes the
metric on $T$. The infimum is always attained.  If $\ell(\phi)>0$
there is a unique $\phi$-invariant line (=isometric image of $\R$),
called the {\it axis} of $\phi$, and the restriction of $\phi$ to this
line is translation by $\ell(\phi)$. In this case $\phi$ is said to be
{\it hyperbolic}. If $\ell(\phi)=0$, then $\phi$ fixes a non-empty
subtree of $T$ and is said to be {\it elliptic}.

\begin{exercise} Let $\phi$ and $\psi$ be two isometries of an \rt\ $T$.
If they are both elliptic with disjoint fixed point sets, then the
composition $\psi\phi$ is hyperbolic, and $\ell(\psi\phi)$ is equal
to twice the distance between $Fix(\phi)$ and $Fix(\psi)$. If both
$\phi$ and $\psi$ are hyperbolic and their axes are disjoint, then
$\psi\phi$ is hyperbolic, the translation length is equal to the
sum of the translation lengths of $\phi$ and $\psi$ plus twice the
distance between the axes of $\phi$ and $\psi$, and the axis of
$\psi\phi$ intersects both the axis of $\phi$ and of $\psi$.
\end{exercise}

\begin{exercise} \label{two}
If $\{T_i\}_{i\in I}$ is a 
finite collection of subtrees
of an \rt\ $T$ such that all pairwise intersections are non-empty,
then the intersection of the whole collection is non-empty.
\end{exercise}

Now let $G$ be a group acting by isometries on an \rt\ $T$. The action
is {\it non-trivial} if no point of $T$ is fixed by the whole group.
It is {\it minimal} if there is no proper $G$-invariant subtree.

\begin{exercise}
Use Exercise \ref{two} to show that whenever a finitely generated
group acts non-trivially on an \rt, then some elements of the group
are mapped to hyperbolic isometries. Construct a (simplicial)
counterexample to this statement when ``finitely generated'' is
omitted from the hypotheses.
\end{exercise}

\begin{prop} \label{countably many axes}
Assume that $G$ is finitely generated and that the action of $G$
on $T$ is non-trivial. Then $T$ contains a unique $G$-invariant subtree
$T'\subset T$ such that the action restricted to $T'$ is minimal.
Further, $T'$ is the union of at most countably many lines.
\end{prop}

\begin{proof} Let $T'$ be the union of the axes of hyperbolic elements
in $T$. The only fact that needs a proof is that $T'$ is non-empty and
connected, and this follows from the exercises above.
\end{proof}

We will often replace a given \rt\ with the minimal subtree without
saying so explicitly.

\subsection{\dh\ spaces -- a review}

We recall the notion of ``negative curvature'' for metric spaces, due
to M. Gromov \cite{mg:hyperbolic}. For an inspired exposition see
\cite{gh:gromov}.

\begin{definition}
Let $(X,d)$ be a metric space and $*\in X$ a basepoint. For $x,y\in X$ define
$(x\cdot y)=\frac{1}{2}(d(*,x)+d(*,y)-d(x,y))$. For $\delta\geq 0$ we say that
$(X,*,d)$ is {\it $\delta$-hyperbolic} if for all $x,y,z\in X$ we have
$$(x\cdot
y)\geq \min ((x\cdot z),(y\cdot z))-\delta.$$
\end{definition}

\begin{example}
Let $X$ be an $\R$-tree. Then $(x\cdot y)$ equals the distance between
$*$ and the segment $[x,y]$. Further, if $x,y,z\in X$, then $(x\cdot
y)\geq \min ((x\cdot z),(y\cdot z))$. Thus $\R$-trees are
0-hyperbolic spaces. The converse is given in the next lemma.
\end{example}

\begin{example}
Hyperbolic space $\H^n$ and any complete simply-connected
Riemannian manifold with sectional curvature $\leq -\epsilon<0$ is \dh\
for some $\delta=\delta(\epsilon)$.
\end{example}

If $(X,*,d)$ is \dh\ and if $*'$ is another basepoint, then
$(X,*',d)$ is $2\delta$-hyperbolic. It therefore follows that the
notions of ``0-hyperbolic'' and ``hyperbolic'' (i.e. \dh\ for some
$\delta$) don't depend on the choice of the basepoint. 

A finitely generated group $G$ is {\it word-hyperbolic}
\cite{mg:hyperbolic} if the word metric on $G$ with respect to a
finite generating set is hyperbolic. This notion is independent of the
choice of the generating set. 

The classification of isometries of hyperbolic
spaces is more subtle than in the case of trees. Let $\phi:X\to X$ be
an isometry of a hyperbolic metric space. The {\it translation length}
can be defined as the limit
$$\ell(\phi)=\lim_{i\to\infty}\frac{1}{i}\inf_{x\in X}
d(x,\phi^i(x)).$$ For a reasonable classification into hyperbolic,
elliptic, and parabolic isometries it is necessary to assume something
about $X$, e.g. that it is a geodesic metric space (any two points can
be joined by a geodesic segment), or perhaps something weaker that
guarantees that $X$ does not have big holes.  We will only need the
following special case. If $G$ is a word-hyperbolic group and $g\in G$
an element of infinite order, then left translation $t_g:G\to G$ by
$g$ has an {\it axis}, namely the set of points in $G$ moved a
distance $\leq \ell(t_g)+10\delta$. ($G$ is $\delta$-hyperbolic.) This
set is quasi-isometric to the line.

\subsection{``Connecting the dots'' lemma}

\begin{lemma} \label{connecting-the-dots}
Let $(X,*,d)$ be a 0-hyperbolic metric space. Then there exists an $\R$-tree
$(T,d_T)$ and an isometric embedding $i:X\to T$ such that
\begin{enumerate}
\item no proper subtree of $T$ contains $i(X)$, and
\item if $j:X\to T'$ is an isometric embedding of $X$ into an $\R$-tree $T'$,
then there is a unique isometric embedding $k:T\to T'$ such that $ki=j$.
\end{enumerate}
\end{lemma}

In particular, $T$ is unique up to isometry. Further, if a group $G$ acts by
isometries on $X$, then the action extends to an isometric action on $T$.

\begin{proof}
If $i:X\to T$ is an isometric embedding as in (1), then $T$ is the
union of segments of the form $I_x=[i(*),i(x)]$ for $x\in X$. The
length of $I_x$ is equal to $d(*,x)$, and two such segments $I_x$ and
$I_y$ overlap in a segment of length $(x\cdot y)$. This suggests the
construction of $T$. Start with the collection of segments
$I_x=[0,d(*,x)]$ for $x\in X$ and then identify $I_x$ and $I_y$ along
$[0,(x\cdot y)]$. For details, see e.g. \cite{otal:hyperbolization}.
\end{proof}

\section{How do \rt s arise?}

Isometric actions of a group $G$ on an \rt\ arise most often as
the Gromov-Hausdorff limits of a sequence of isometric actions of 
$G$ on a negatively curved space $X$. The construction is due independently
to F. Paulin \cite{fp:trees} and to M. Bestvina \cite{mb:degenerations}.
See also the expository article \cite{bridson-swarup:paulin}.

We will make the formal
definition in terms of the projectivized space of
equivariant pseudometrics on $G$.

\subsection{Convergence of based $G$-spaces}
Let $G$ be a discrete group. By a {\it $G$-space} we mean a pair
$(X,\rho)$ where $X$ is a metric space and $\rho:G\to Isom(X)$ is a
homomorphism (an {\it action})
to the group of isometries of
$X$. A {\it based $G$-space} is a triple $(X,*,\rho)$ where
$(X,\rho)$ is a $G$-space and $*$ is a basepoint in $X$ that is not
fixed by every element of $G$.

Recall that a {\it pseudometric} on $G$ is a function $d:G\times G\to
[0,\infty)$ that is symmetric, vanishes on the diagonal, and satisfies
the triangle inequality. Let $\cal D$ denote the space of all
pseudometrics (``distance functions'') on $G$ that are not identically
0, equipped with compact-open topology. We let $G$ act on $G\times G$
diagonally, and on $[0,\infty)$ trivially, and consider the subspace
${\cal ED}\subset \cal D$ of $G$-equivariant pseudometrics. Scaling
induces a free action of $\R^+$ on $\cal ED$, and we denote by $\cal
PED$ the quotient space, i.e. the space of {\it projectivized
equivariant distance functions} on $G$. A pseudometric on $G$ is
$\delta$-hyperbolic if the associated metric space is
$\delta$-hyperbolic (the class of the identity element is taken to be
the basepoint). 

A based $G$-space $(X,*,\rho)$ induces an equivariant pseudometric
$d=d_{(X,*,\rho)}$ on
$G$ by setting
$$d(g,h)=d_X(\rho(g)(*),\rho(h)(*))$$
where $d_X$ denotes the distance function in $X$.
If the stabilizer under $\rho$ of $*$ is trivial, then $G$ can be
identified with the orbit of $*$ via $g\leftrightarrow \rho(g)(*)$,
and $d_{(X,*,\rho)}$ is the distance induced by $d_X$. We work with
pseudometrics to allow for the possibility that distinct elements of
$G$ correspond to the same point of $X$.

\begin{definition}
We say that a sequence $(X_i,*_i,\rho_i)$, $i=1,2,3,\cdots$ of
based $G$-spaces {\it converges} to the based $G$-space $(X,*,\rho)$
and write
$$\lim_{i\to\infty}(X_i,*_i,\rho_i)=(X,*,\rho)$$
provided $[d_{(X_i,*_i,\rho_i)}]\to [d_{(X,*,\rho)}]$ in $\cal PED$.
\end{definition}

\subsection{Example: Flat tori}
To illustrate this, let us take $G=\Z \times \Z$ and $X=E^2$, the
Euclidean plane. We will obtain actions of $G$ on
the real line as limits of discrete actions of $\Z\times\Z$ on $E^2$.
The group
$\Z\times\Z$ can act on
$E^2$ by isometries in many different ways. We will only consider
discrete isometric actions, and those consist necessarily of
translations and form the universal covering group of a flat 2-torus.
Two such actions of
$\Z\times\Z$ will be considered equivalent if there is a similarity of
$E^2$ conjugating one action to the other, i.e. if the corresponding
(marked) tori are conformally equivalent. It is convenient to
identify the group of translations of $E^2$ with $\C$. Thus two
actions $\rho_1,\rho_2:G\to \C$ are equivalent if there is a
complex number $\alpha$ such that $\rho_2(g)=\alpha \rho_1(g)$
for
all $g\in \Z\times\Z$ or $\rho_2(g)=\alpha \overline{\rho_1(g)}$ for
all
$g\in \Z\times\Z$. Each equivalence class $[\rho]$ is uniquely
determined by the complex-conjugate pair $\{ z,\overline z\}$ where
$z=\frac{\rho(0,1)}{\rho(1,0)}$, and
thus the set of all equivalence classes can be identified with
the upper half-plane $\{ z\in \C| Im(z)>0\}$. 

Let us take a sequence $\{ z_n\}$ of points in the upper
half-plane and let $\rho_n$ be a representative of the equivalence
class determined by $z_n$. Fix a basepoint $*$ in
$E^2$.

\begin{prop} Suppose that $z_n\to r\in \R\cup \{ \infty\}$. 
Then the sequence $(E^2,*,\rho_n)$ converges to the (unique up to
scale) based $G$-space $(\R,0,\rho)$ where the action $\rho$ consists
of translations and $\frac{\rho(0,1)}{\rho(1,0)}=r$.
\end{prop}

\begin{proof}
Suppose for concreteness that $r\in\R$. We take $\rho_n$ so that
$\rho_n(1,0)=1$ and $\rho_n(0,1)=z_n$. Then $\rho_n(g)\to \rho(g)\in
\R\subset\C$ and the claim follows.
\end{proof}

Thus this
construction recovers the usual compactification of the upper
half-plane by the circle $\R\cup \{ \infty\}$.

If in this example we replace $E^2$ by $E^n$ and $\Z \times \Z$ by
$\Z^n$, the same construction would produce actions of $\Z^n$ by
translations on $E^m$, $0<m<n$ and would provide an equivariant
compactification of the symmetric space $SL_n(\R)/SO_n(\R)$.

\subsection{\rt s as limits of based \dh\ $G$-spaces}
The main reason for interest in \rt s is the following result. Note
that if a sequence of based $G$-spaces converges, the limit is far
from being unique. In special situations we can take the limit to
be an \rt.

\begin{thm} \label{compactness}
Let $(X_i,*_i,\rho_i)$ be a convergent sequence of
based $G$-spaces. Assume that
\begin{enumerate}
\item there exists $\delta\geq 0$ such that every $X_i$ is \dh, and
\item there exists $g\in G$ such that the sequence
$d_{X_i}(*,\rho_i(g)(*))$ is unbounded.
\end{enumerate}
Then there is a based $G$-tree $(T,*)$ and an isometric action
$\rho:G\to Isom(T)$ such that $(X_i,*_i,\rho_i)\to (T,*,\rho)$.
\end{thm}

\begin{proof}
The limiting pseudo-metric $d$ on $G$ is 0-hyperbolic, as it is the
limit of pseudo-metrics $\frac{d_{(X_i,*_i,\rho_i)}}{d_i}$ with
$d_i\to\infty$ (by (2)) and the $i$th pseudo-metric is
$\frac{\delta}{d_i}$-hyperbolic (by (1)). Now apply the
connecting-the-dots Lemma \ref{connecting-the-dots} to the induced
metric space.
\end{proof}

\begin{exercise} \label{axes off}
Let $G=\Z$, $X_i=\H^2$ (hyperbolic plane), and the
representation $\rho_i$ sends the generator $1\in \Z$ to a
hyperbolic isometry whose translation length is $1$ and whose axis
passes at distance $i$ from the basepoint in $\H^2$.
Show that the limiting \rt\ $T$ can be identified with the cone on $\Z$
and the limiting action $\rho$ is the translation action on $\Z$ coned
off. The basepoint is a point in $\Z$.
\end{exercise}

\begin{example} \label{tracks}
Let $f:F_n\to F_n$ be an automorphism of the free group
$G=F_n=<x_1,\cdots,x_n>$ of rank $n$ that sends each basis element to
a ``positive word'', i.e. a product of basis elements (not involving
their inverses). Suppose that $\lambda>0$ is the unique eigenvalue of
the abelianization of $f$, viewed as an automorphism of $\Z^n$, with a
corresponding eigenvector with non-negative coordinates
$a_1,\cdots,a_n$. For $X_i$ take $F_n$ with the word metric, the
basepoint is 1, and let $\rho_i$ be the representation that sends
$g\in F_n$ to the left translation by $f^i(g)$. The
scaling factor can be taken to be $d_i=\lambda^i$. If $\lambda>1$ we
obtain in the limit an action $\rho$ of $F_n$ on an \rt. The
positivity requirement was imposed to ensure that for some $g$ the
sequence of lengths of $f^i(g)$ grows at the ``top speed'', i.e. as
$(const)\lambda^i$. The limiting tree can be described quite
explicitly. For example, if $g$ is a positive word, then the distance
between the basepoint $*\in T$ and its image under $\rho(g)$ is
$k_1a_1+\cdots+k_na_n$, where $k_j$ is the number of times the
generator $x_j$ appears in the word $g$.

More generally, this construction can be performed with the
train-track maps of \cite{bh:tracks}. With the right choice of a
train-track map one obtains free nonsimplicial isometric actions of
the free group $F_n$ ($n>2$) on \rt s (see \cite{ps:dendrology2}).
\end{example}

\subsection{Finding approximate subtrees}
We now assume that we are in the situation of Theorem
\ref{compactness} and we examine the limiting tree in more
detail. Thus we assume that there is a sequence 
$d_i\to\infty$ such that
$$d_T(\rho(g)(*),\rho(h)(*))=\lim_{i\to\infty}\frac{d_{X_i}
(\rho_i(g)(*_i),\rho_i(h)(*_i))}{d_i}.$$ 
If $x$ is a point in $T$ that belongs to the orbit of $*$, then $x$
can be ``approximated'' by the corresponding point $x_i\in X_i$ in the
orbit of $*_i$. If two points in the orbit of $*$ coincide, then the
corresponding points in $X_i$ are ``close'' (more precisely, the
distance between them divided by $d_i$ goes to 0 as $i\to\infty$).

We now extend this discussion to all points of $T$. We will assume in
addition that each $X_i$ is a {\it geodesic metric space}, i.e. that
any two points $x,y\in X_i$ are joined by a geodesic segment.

Let $x\in T$ be an arbitrary point. Fix a
finite subset $F\subset G$ and define a set $X_i(F,x)\subset X_i$
``approximating'' $x$ to be the set of all points $x_i\in X_i$ that
can be constructed as follows. Choose $g,h\in F$ so that $x$ is on the
segment in $T$ connecting $\rho(g)(*)$ and $\rho(h)(*)$ and choose a
geodesic segment in $X_i$ connecting $\rho_i(g)(*_i)$ and
$\rho_i(h)(*_i)$, and let $x_i$ be the point on this segment that
divides it in the same ratio as the point $x$ divides the segment
$[\rho(g)(*),\rho(h)(*)]$.

Of course, it might happen that $X_i(F,x)=\emptyset$ if there are no
$g,h\in F$ as above. The following proposition summarizes
the basic properties of this construction:

\begin{prop} \label{approximate trees}
\hfill
\begin{enumerate}
\item {\it equivariance:} $X_i(gF,\rho(g)x)=\rho_i(g)X_i(F,x)$.
\item {\it monotonicity:} If $F\subset F'$ then $X_i(F,x)\subset
X_i(F',x)$.
\item {\it small diameter:} $\frac{1}{d_i}diam\ X_i(F,x)\to 0$ as
$i\to\infty$.
\item {\it metric convergence:} Let $x,y\in T$. Then for all finite 
$F\subset G$ and all choices
$x_i\in X_i(F,x)$, $y_i\in X_i(F,y)$ we have
$$\frac{1}{d_i}d_{X_i}(x_i,y_i)\to d_T(x,y).$$
\item {\it non-triviality:} For every $x\in T$ there is a 2-element
set $F\subset G$ such that $X_i(F,x)\neq\emptyset$ for all $i$
\end{enumerate}
\end{prop}

\begin{proof}
Items 1 and 2 follow directly from the definition. Item 3 is an
exercise in $\delta$-hyperbolic geometry. Item 4 also follows directly
from definitions if $F=\{ g,h\}$ so that the segment in $T$ joining
$\rho(g)(*)$ and $\rho(h)(*)$ contains both $x$ and $y$ (the existence
of such $g,h$ follows from item 1 of Lemma \ref{connecting-the-dots},
which also implies item 5). The general case then follows from 2 and 3
by enlarging $F$.
\end{proof}

\subsection{Selecting the basepoint and the Compactness Theorem}
We now assume that an action $\rho:G\to Isom(X)$ is given, and we
consider the problem of locating a ``most centrally located point''
for this action. We will then use this point as the basepoint. In
Example \ref{axes off} the basepoint $*_i$ should be chosen on the
axis of $\rho_i(1)$, and then the limiting action would be
nontrivial. 

The problem of finding a good basepoint has a satisfactory solution
when $G$ is finitely generated and $X$ is a {\it proper} \dh\ metric
space, and this is what we assume from now on. (A metric space is {\it
proper} if closed metric balls are compact.) We also fix a finite
generating set $S\subset G$. Let $F=F_{S,\rho}:X\to [0,\infty)$ be the
function defined by
$$ F(x)=\max_{g\in S}d_X(x,\rho(g)(x)).$$
The following lemma is an exercise.

\begin{lemma} Assume that $\rho:G\to Isom(X)$ is non-elementary (i.e.
it does not fix a point at infinity). Then
$F:X\to [0,\infty)$ is a proper map. In particular, $F$ attains its
global minimum.\qed
\end{lemma}

We call a point $x\in X$ {\it centrally located (with respect to the
action $\rho:G\to Isom(X)$ and the generating set $S$)} if
$F$ attains its global minimum at $x$.

\begin{prop}
Suppose that under the hypotheses of Theorem \ref{compactness} each
$X_i$ is proper and that the basepoints $*_i$ are centrally located
(with respect to $\rho_i$ and a fixed finite generating set $S$ for
$G$). Then the limiting action $\rho:G\to Isom(T)$ does not have global
fixed points.
\end{prop}

\begin{proof}
We can take $d_i=\max_{g\in S}d_{X_i}(*_i,\rho_i(g)(*_i))$. Suppose
$x\in T$ is a global fixed point. Choose a finite subset $F\subset G$
so that $X_i(F,x)\neq\emptyset$. We will argue that for any $x_i\in
X_i(F,x)$ and any $g\in S$ we have
$\frac{1}{d_i}d_{X_i}(x_i,\rho_i(g)(x_i))\to 0$ as $i\to \infty$,
contradicting (for large $i$) the assumption that $*_i$ is centrally
located. Indeed, $\rho(g)(x)=x$ coupled with equivariance property
implies $X_i(gF,x)=\rho_i(g)X_i(F,x)$, so by monotonicity both $x_i$
and
$\rho_i(g)(x_i)$ belong to $X_i(gF\cup F,x)$, so the claim follows
from the small diameter property.
\end{proof}

\begin{thm}[Compactness Theorem]
Suppose that $(X,d)$ is a proper \dh\ metric space and that $\rho_i:G\to
Isom(X)$ a sequence of non-elementary representations of a finitely
generated group $G$. Assume that the group $Isom(X)$ acts cocompactly
on $X$, i.e. that there is a compact subset $K\subset X$ whose
$Isom(X)$-translates cover $X$. Then one of the following holds,
possibly after passing to a subsequence.
\begin{enumerate}
\item There exist isometries $\phi_i\in Isom(X)$ such that the sequence
of conjugates $\rho_i^{\phi_i}$ converges in the compact-open topology
to a representation $\rho:G\to Isom(X)$.
\item For each $i$ there exists a centrally located point $x_i\in X$
for the representation $\rho_i$ such that the sequence of based
$G$-spaces $(X,x_i,\rho_i)$ converges to an action of $G$ on an \rt\
$T$ without global fixed points.\end{enumerate}\end{thm}

\begin{proof}
Let $x_i$ be a centrally located point for $\rho_i$. If the sequence
$d_i=\max_{g\in S}d_{X}(x_i,\rho_i(g)(x_i))$ converges to infinity,
then item 2 holds, by the preceding proposition. Otherwise, after
passing to a subsequence, the $d_i$'s are uniformly bounded. In that
case choose $\phi_i\in Isom(X)$ that sends $x_i$ into $K\subset X$ and
apply Arzela-Ascoli to the conjugates $\rho_i^{\phi_i}$ to see that
item 1 holds in this case.
\end{proof}

In the situation (2) of the Compactness Theorem, assuming that $X$ is a
geodesic metric space or a hyperbolic group $G$ equipped with a word
metric, it can also be argued that the translation length
$\ell(\rho(g))$ is equal to the limit
$\lim_{i\to\infty}\frac{\ell(\rho_i(g))}{d_i}$.

\subsection{Arc stabilizers}
We now investigate, in the situation (2) of the Compactness Theorem,
the arc stabilizers in the limiting action $\rho:G\to Isom(T)$. We
restrict ourselves to two frequently encountered settings, when the
arc stabilizers turn out to be ``elementary''.

Many reasonable groups, such as linear groups, satisfy the so called
``Tits Alternative''. This means that their subgroups are either
``small'' (virtually solvable) or ``large'' (contain a nonabelian free
group).  Word-hyperbolic groups satisfy a strong form of the Tits
alternative: any subgroup either contains a nonabelian free group or
it is virtually cyclic (elementary). Accordingly, an action of a group
on an $\R$-tree is said to be {\it small} if it is non-trivial (there
are no global fixed points), minimal, and all arc stabilizers are
small.

\begin{prop} \label{small}
Let $H\subset G$ be the stabilizer under $\rho$ of a non-degenerate arc
in
$T$.
\begin{enumerate}
\item If each $X_i$ is a copy of the Cayley graph $\Gamma$ of a word
hyperbolic group with respect to a fixed finite generating
set, and each $\rho_i:G\to Isom(\Gamma)$ is a free action whose image
is contained in the subgroup consisting of left translations, then $H$
is virtually cyclic.
\item If each $X_i$ is a copy of a fixed rank 1 symmetric space $\H$
(real, complex, quaternionic hyperbolic space, or the Cayley plane),
and each $\rho_i$ is discrete and faithful, then $H$ is virtually
nilpotent.
\end{enumerate}
\end{prop}

\begin{proof}
Let $[a,b]\subset T$ be a non-degenerate
segment fixed by $H$ (under the action by
$\rho$). Choose a sufficiently large finite subset
$F\subset G$ and points
$a_i\in X_i(F,a)$ and
$b_i\in X_i(F,b)$. Let $c_i$ be the midpoint on a geodesic segment
$\sigma_i$ connecting $a_i$ and $b_i$. 

(1) Say $\Gamma$ is \dh. The key claim is that if $h,k\in H$ then, for
large $i$, the left translation $\rho_i([h,k])$ moves $c_i$ to a point
at distance $<20\delta$ from $c_i$. There is an upper bound to the
number of left translations of $\Gamma$ that move a given point a
distance $\leq 20\delta$. Since the commutators $[h,k]$ for $h,k\in H$
generate the commutator subgroup $[H,H]$ of $H$, it follows from the
freeness assumption that $[H,H]$ is finitely generated and in
particular $H$ is not a nonabelian free group. Since the same
argument can be applied to any subgroup of $H$, we conclude that $H$
does not contain a nonabelian free group, and hence it is virtually
cyclic. 

The idea of proof of the above key claim is that $\rho_i(h)$ and
$\rho_i(k)$ map $\sigma_i$ to a geodesic segment whose endpoints are
within $\frac{1}{100}length(\sigma_i)$ of the endpoints of $\sigma_i$,
and so these segments, except near the endpoints, run within $2\delta$
of $\sigma_i$, i.e. $\rho_i(h)$ and $\rho_i(k)$ can be thought of
(modulo small error) as translating along $\sigma_i$.  Consequently,
the commutator $\rho_i([h,k])$ fixes $\sigma_i$ (modulo small error
and away from the endpoints). Details are in \cite{mb:degenerations}
and \cite{fp:trees}.

(2) The proof here is a modification of (1), plus the Margulis lemma.
Let $\mu$ be the Margulis constant for $\H$, so that if a discrete
group of isometries of $\H$ is generated by isometries that move a
point $x_0\in \H$ a distance $<\mu$, then the group is virtually
nilpotent. Arguing as in the key claim above, one can show that if
$h,k\in H$, then for large $i$ the isometry $\rho_i([h,k])$ moves
$c_i$ a distance $<\mu$. It then follows that every finitely
generated subgroup of $[H,H]$ is virtually nilpotent and so
$\rho_i(H)$ must be elementary (i.e. virtually nilpotent).
\end{proof}

\subsection{Stable actions}

\begin{definition}
Suppose a group $G$ is acting isometrically on an \rt\ $T$. A subtree
of $T$ is {\it non-degenerate} if it contains more than one point. A
non-degenerate subtree $T_1\subset T$ is said to be {\it stable} (with
respect to the action) if for every non-degenerate subtree $T_2\subset
T_1$ we have the equality $Fix(T_1)=Fix(T_2)$ of pointwise
stabilizers. The group action on $T$ is {\it stable} if it is non-trivial,
minimal, and every
non-degenerate tree in $T$ contains a stable subtree.
\end{definition}

Group actions that tend to arise in practice are stable. For example,
small actions of hyperbolic groups are stable. More generally, if the
collection of arc stabilizers satisfies the ascending chain condition,
then the (non-trivial and minimal) group action is stable. Note that
if two stable subtrees of $T$ have a non-degenerate intersection, then
their union is a stable subtree. In particular, each stable subtree is
contained in a unique maximal stable subtree.

The study of stable actions quickly reduces to the study of actions
with trivial arc stabilizers (see Corollary 5.9 of \cite{bf:stable}).
To see the idea, assume that $T$ is covered by
maximal stable subtrees $\{T_i\}_{i\in I}$. Note that $T_i\cap T_j$ is at
most a point for $i\neq j$.  Now construct a simplicial tree $S$ as
follows. There are two kinds of vertices in $S$. There is a vertex for
each maximal stable subtree $T_i$, and there is a vertex for each
point of $T$ that equals the intersection of distinct maximal stable
subtrees. An edge is drawn from a vertex $v$ of the first kind,
determined by $T_i$, to the vertex $w$ of the second kind, determined
by $x\in T$, precisely when $x\in T_i$. The group $G$ acts
simplicially, without inversions of edges, on $S$. The stabilizer
$Fix_S([v,w])$ of the edge $[v,w]$ described above 
fixes a point of $T$
and the underlying assumption is
that we understand arc and point stabilizers in $T$. 
We then appeal to Bass-Serre theory
\cite{se:trees} to conclude that either $G$ splits over an edge
stabilizer in $S$ or that $G$ fixes a vertex of $S$. In the latter
case, in view of nontriviality and minimality of the action of $G$ on $T$,
it follows that $T$ itself is a stable tree, so after factoring out the
kernel of the action, the induced action has trivial arc stabilizers.

\section{Measured laminations on 2-complexes}

We now review the basics of measured laminations. For more information
and details the reader is referred to \cite{ms:valuations2}.

\begin{definition} A closed subset $\Lambda$ of a locally path-connected
metrizable space $X$
is a {\it lamination} if every point $x\in\Lambda$ has a neighborhood $U$
such that the pair $(U,U\cap \Lambda)$ is homeomorphic to the pair
$(V\times (0,1),V\times C)$ for some topological space $V$ and some compact
totally disconnected subset $C\subset (0,1)$. Such a homeomorphism is
called a {\it chart}. The path components of
$\Lambda$ are called {\it leaves}.\end{definition}

If $X$ is a closed manifold, any codimension 1 submanifold is a
lamination. More typically, the set $C$ in the definition is the Cantor set.

\begin{example} Let $X$ be a closed hyperbolic surface. Let $\gamma_i$ be
a sequence of simple closed geodesics in $X$.  After possibly passing
to a subsequence, this sequence converges in the Hausdorff metric to a
closed subset $\Lambda$ of $X$. One can check \cite{cb:surfaces}
that $\Lambda$ is a lamination, and that the leaves of $\Lambda$ are
simple geodesics (closed or biinfinite).  Such $\Lambda$ is called a
{\it geodesic lamination}.\end{example}

\begin{definition} Let $\Lambda\subset X$ be a lamination and $\alpha:[a,b]\to
X$ a path in $X$ such that $\alpha(a),\alpha(b)\notin\Lambda$. We say
that $\alpha$ is {\it transverse} to $\Lambda$ if for every $t\in
[a,b]$ with $\alpha(t)\in \Lambda$ there is a chart
$h:(U,U\cap\Lambda)\to (V\times (0,1),V\times C)$ at $\alpha(t)\in X$
such that the map $pr_{(0,1)}h\alpha$ is a
local homeomorphism at $t$.  \end{definition}

\begin{definition} A {\it transverse measure} on a lamination
$\Lambda\subset X$ is a function $\mu$ that assigns a nonnegative real number
$\mu(\alpha)$ to every path $\alpha$ transverse to $\Lambda$ and satisfies
the following properties.
\begin{enumerate}
\item If $\alpha$ is the concatenation of paths $\beta$ and $\gamma$
both of which are transverse to $\Lambda$, then
$\mu(\alpha)=\mu(\beta)+\mu(\gamma)$.
\item Every $x\in\Lambda$ has a chart $(U,U\cap\Lambda)\approx (V\times
(0,1),V\times C)$ and there is a Borel measure $\nu$ on $(0,1)$
supported on $C$ such that for any path $\alpha:[a,b]\to U$ with
endpoints outside $\Lambda$ that projects 1-1 to an interval in
$(0,1)$, the measure $\mu(\alpha)$ equals the $\nu$-measure of the projection.
\end{enumerate}

The number $\mu(\alpha)$ is the {\it measure} of $\alpha$.


\begin{exercise} \label{homotopic}
If two paths are homotopic
through paths transverse to $\Lambda$, then they have the same measure.
If a path is reparametrized, its measure does not change.
\end{exercise}

The {\it support} of $\mu$ is the complement of the set of points such
that $\mu(\alpha)=0$ whenever the image of $\alpha$ is contained in a
sufficiently small neighborhood of the point.  A lamination is {\it
measured} if it is equipped with a transverse measure.
\end{definition}

The support of $\mu$ is always a sublamination of $\Lambda$. We say that
$\Lambda$ has {\it full support} if the support is all of $\Lambda$.

\begin{example}\label{hypsurf} 
Let $\Lambda$ be a geodesic lamination on a hyperbolic
surface $X$. If $\Lambda$ is the finite union of simple closed curves,
any transverse measure assigns a nonnegative real number, the {\it
multiplicity} to each leaf, and the measure of any path transverse to
$\Lambda$ is the geometric intersection number with $\Lambda$, counted
with multiplicity.
Conversely, any such assignment determines a transverse measure.
Now suppose that $\ell$ is an infinite leaf of $\Lambda$. We will
construct a transverse measure on $\Lambda$, called the {\it counting}
or the {\it hitting} measure. Triangulate the surface $X$ so that the
vertices are in the complement of $\Lambda$, all edges are geodesic
segments,
and each triangle is contained in a chart for $\Lambda$. For each edge
$e$ the intersection $e\cap\Lambda$ is totally disconnected.  Choose a
point in each component of $Int\ e\setminus \Lambda$. A transverse
measure on $\Lambda$ is determined by its values on the subintervals
of the edges $e$ with endpoints in the selected countable set.
Conversely, if $\mu$ is defined on these countably many special
intervals and the following two conditions hold, then $\mu$ extends
uniquely to a transverse measure on $\Lambda$:
\begin{enumerate}
\item{(additivity)} If a special interval $I$ is the concatenation of two
special subintervals $I_1$ and $I_2$, then $\mu(I)=\mu(I_1)+\mu(I_2)$.
\item{(compatibility)} If $I_1$ and $I_2$ are special intervals belonging
to two edges of the same triangle $T$ in the triangulation, and if there
is an embedded quadrilateral in $T$ with two opposite sides $I_1$ and
$I_2$, and the other two opposite sides disjoint from $\Lambda$, then
$\mu(I_1)=\mu(I_2)$.
\end{enumerate}
We will now construct a hitting measure on $\Lambda$.
Choose a sequence of longer and longer closed subintervals
$L_1,L_2,\cdots$ of the leaf $\ell$. For a special interval $I$ define
$$\mu_i(I)=\frac{N(L_i,I)}{N(L_i,X^{(1)})}$$ where $N(L_i,I)$ is the
number of intersection points in $L_i\cap I$ and similarly
$N_i=N(L_i,X^{(1)})$ is the number of intersection points between $L_i$
and the 1-skeleton. Since $\ell$ is an infinite leaf, we have $N_i\to\infty$.
Additivity and compatibility hold approximately for $\mu_i$, i.e. in both
cases the difference between the left-hand and the right-hand side is in
the interval $[-\frac{1}{N_i},\frac{1}{N_i}]$. Using a diagonalization process,
pass to a subsequence if necessary so that $\lim \mu_i(I)$ exists for
each special interval $I$, and set the limiting value equal to $\mu(I)$.
The support of $\mu$ is generally smaller than $\Lambda$.
\end{example}

\begin{example}\label{torus}
Let $f:\R^n\to\R$ be a linear map that is injective when restricted to
$\Z^n$ and consider the foliation of $\R^n$ by the level sets of $f$
and the induced foliation $\cal F$ on the torus $T^n=\R^n/\Z^n$. The map $f$
also defines a transverse measure on the two foliations. There is a
standard way of converting the measured foliation $\cal F$ to a measured
lamination $\Lambda$ (and vice-versa); indeed, \cite{bf:stable} is written in
the language of foliations. More precisely, there is a map $T^n\to T^n$
whose point-preimages are arcs and points and the preimage of each leaf of
$\cal F$ is either a leaf of $\Lambda$ or the closure of a complementary
component of $\Lambda$. This map is modeled on the Cantor function $[0,1]\to
[0,1]$, which converts the foliation of $[0,1]$ by points to the lamination
on $[0,1]$ whose underlying set is the Cantor set.
\end{example}

\subsection{Sacksteder's Theorem}

We say that two paths $\gamma$ and $\delta$ transverse to a lamination
$\Lambda$ are {\it pushing equivalent} if, after possibly
reparametrizing one, they are homotopic through paths transverse to
$\Lambda$. If $\Lambda$ is equipped with a transverse measure $\mu$,
then by Exercise \ref{homotopic} we have $\mu(\gamma)=\mu(\delta)$. In
particular, a measured lamination of full support satisfies the
following {\it non-nesting condition}:
\vskip .1cm
\noindent {\it If $\gamma:[a,b]\to X$ and $\delta:[c,d]\to X$ are
pushing-equivalent and
$\gamma$ is a subpath of $\delta$ (i.e. $\gamma=\delta|[a,b]$), then
$\delta([c,d]\setminus [a,b])\subset X\setminus \Lambda$.}
\vskip .1cm
There is a remarkable converse, due to R. Sacksteder.

\begin{thm} \cite{sack:sack} \label{sack}
Suppose $X$ is compact and $\Lambda\subset X$ is a lamination on $X$
satisfying the above non-nesting condition. Then there is a non-trivial
transverse measure on $\Lambda$, possibly not of full support.\end{thm}

It is easy to construct examples of non-nesting laminations on compact spaces
that do not support a transverse measure. For example, take a geodesic
measured lamination on a hyperbolic surface and replace a noncompact leaf
with a parallel family of leaves. For an \rt\ version of Sacksteder's
theorem, see 
\cite{levitt:sack}.

\subsection{Decomposition into minimal and simplicial components}

We say that a lamination $\Lambda\subset X$ is {\it simplicial} if there is a
leaf
$\ell$ of $\Lambda$, a closed neighborhood $N$ of $\Lambda$ in $X$ and a map
$N\to\ell$ which is an $I$-bundle and whose restriction to $\Lambda$ is a
bundle map with 0-dimensional fibers. A lamination $\Lambda$ is {\it minimal}
if every leaf of $\Lambda$ is dense in $\Lambda$. When the underlying space
$X$ is compact, the lamination that supports a transverse measure always
decomposes into simplicial and minimal sub-laminations.

\begin{thm} (Theorem 3.2 in \cite{ms:valuations2}) 
Let $X$ be compact and $\Lambda\subset X$ a
lamination that admits a transverse measure with full support. Then
$\Lambda$ is the disjoint union
$\Lambda_1\sqcup\Lambda_2\sqcup\cdots\sqcup\Lambda_n$ with each
$\Lambda_i$ either simplicial or minimal.\end{thm}

On a closed hyperbolic surface, imagine a lamination consisting of two closed
geodesics and a biinfinite geodesic that spirals towards the closed
geodesics, one in each direction. Such a lamination does not decompose into
simplicial and minimal sub-laminations. It is also not hard to show directly
that this lamination does not support a transverse measure; indeed, this
lamination is not even non-nesting.

\subsection{Resolutions}\label{resolutions}

Let $G$ be a finitely presented group, and assume that $G$ is acting
non-trivially and minimally on an \rt\ $T$ (as usual, by
isometries). Since $G$ is finitely presented, there is a finite
simplicial complex $K$ of dimension $\leq 2$ whose fundamental group
is $G$. We now use $T$ to construct a measured lamination $\Lambda$ on
$K$ and an equivariant map $f:\tilde K\to T$ from the universal cover
of $K$ to $T$ that sends leaves of the preimage lamination
$\tilde\Lambda\subset \tilde K$ to points. We refer to this map as a
{\it resolution}. In the case of simplicial trees this construction
has been extensively used by M. Dunwoody (the leaves in this case are
Dunwoody's ``tracks'').

To construct $\Lambda$ and $f$, first choose a countable equivariant
dense subset $D\subset T$ that includes all branch points of $T$
($v\in T$ is a branch point if the tripod, i.e. the cone on 3 points,
can be embedded in $T$ with the cone point mapped to $v$) and that
intersects each arc in a dense set. This is possible by Proposition
\ref{countably many axes}. Then define $f$ on the vertices of $\tilde
K$ so that the map is equivariant and sends each vertex into
$D$. Next, extend $f$ equivariantly to the edges of $\tilde K$. If the
endpoints of a given edge $e$ map to the same point under $f$, then
define $f$ on $e$ to be the constant map.  Otherwise, $f|e$ is chosen
so that it is the Cantor function onto the arc whose boundary is
$f(\partial e)$ with the preimage of each point in $D\cap f(e)$ an arc
and the preimage of every other point in $f(e)$ a single point. The
Cantor set of points in $e$ that don't belong to the interior of a
preimage arc is going to be the set $\tilde\Lambda\cap e$. Finally,
extend $f$ equivariantly to each 2-simplex $\sigma$ of $K$ so that for
each $y\in D\cap f(\sigma)$ the preimage $f^{-1}(y)\cap\sigma$ is a
convex triangle, quadrilateral, or a hexagon with vertices in
$\partial\sigma$ and the preimage of every other point in $f(\sigma)$
is a straight line segment joining two distinct sides of
$\sigma$. These line segments are the components of $\tilde\Lambda\cap
\sigma$. The transverse measure is defined by the requirement that if
$\alpha$ is a path in $\tilde K$ that is transverse to $\tilde\Lambda$
and intersects each leaf at most once, then the measure of $\alpha$ is
the distance in $T$ between the $f$-images of the endpoints of
$\alpha$. This transverse measure is equivariant and descends to a
transverse measure on the induced lamination $\Lambda\subset K$.

\subsection{Dual trees}\label{dual}

There is a construction that to a measured lamination $\Lambda$ on a
finite complex $K$ assigns an \rt\ on which the fundamental group of
the complex acts. Let $\tilde K$ be the universal cover of $K$ and
$\tilde \Lambda$ the induced lamination on $\tilde K$. Define a
pseudometric $d:\tilde K\times\tilde K\to [0,\infty)$ by taking
$$d(x,y)=\inf_{\alpha}\tilde\mu(\alpha)$$ where $\tilde \mu$ is the
induced transverse measure and the infimum runs over all paths that
are transverse to $\tilde\Lambda$ and join $x$ to $y$. It is not
difficult to show that the associated metric space $T$ is an \rt,
called the {\it dual tree}, and that the deck group induces an
isometric action of $G$ on $T$. There is also the natural quotient map
$f:\tilde K\to T$; it is equivariant and maps each leaf and each
complementary component of $\tilde\Lambda$ to a point. 

In general, many different leaves will map to the same point by $f$. For
example, start with a geodesic lamination on a 4 times punctured
sphere and then fill in the punctures. The dual tree in this case is a
single point. It is reasonable to impose the condition that $f$
restricted to each edge $e$ of $\tilde K$ is the Cantor function that
collapses precisely the closures of complementary components of
$e\cap\Lambda$ in $e$. This condition is automatically satisfied when
$\Lambda$ arises as in the construction of a resolution. For the lack
of a better term, we say that $f$ is {\it locally injective} if it
satisfies this condition.

\vskip .2cm
\noindent {\it Questions.} Assume that $f$ is locally injective.  Is
the infimum above always realized by a ``minimizing'' path $\alpha$?
Can $f$ map distinct leaves of $\tilde \Lambda$ that do not belong to
the closure of the same complementary component to the same point?
\vskip .2cm

\begin{example} Let $\Lambda$ be a geodesic measured lamination on a
closed hyperbolic surface such that the measure has full support and
such that the complementary components are simply-connected (see e.g.
\cite{cb:surfaces}). The action of the fundamental group of the
surface on the dual tree is free.\end{example}

Notice that the constructions of a resolution and of the dual tree are
generally not inverses of each other. For example, a free group admits
many interesting non-simplicial actions on \rt s (e.g. via the
construction as in the preceding example applied to a punctured
surface), while the dual of any resolution that uses a bouquet of
circles for $K$ is simplicial.

A resolution $f:\tilde K\to T$ is {\it exact} (see \cite{bf:stable}) if all
point preimages are connected. This is equivalent to the statement that
each point preimage is either a leaf or the closure of a complementary
component. A group action on an \rt\ $T$ is {\it geometric} if it admits an
exact resolution.

Frequently, one encounters the following situation: $f:\tilde K\to T$
is a resolution, $\tilde\Lambda$ the associated lamination, and
$f':\tilde K\to T'$ is the equivariant map to the tree dual to $\tilde
\Lambda$. By construction, we have a factorization $$f=\pi f'$$ for an
equivariant map $\pi:T'\to T$. As remarked above, this map may not be
an isometry. If $f$ is an exact resolution, then $\pi$ is an isometry.
\vskip .2cm
\noindent {\it Question.} If $\pi$ is an isometry, is $f$ an exact
resolution?
\vskip .2cm

It is a consequence of the Rips
machine that if the action on $T$ is stable, then $f':\tilde K\to T'$ is
an exact resolution, so the potential pathologies in the questions above 
don't arise in the stable case.

In general, one can say that if $f$ is not exact, then either $\pi$ is
not an isometry, or there exist two leaves of $\tilde \Lambda$ that
can be joined by a path with arbitrarily small measure, but cannot be
joined by a path of measure 0. In either case, there are two leaves of
$\tilde \Lambda$ such that any path joining them has measure strictly
larger than the distance between their $f$-images $x$ and $y$. One can
then construct a ``better resolution'' as follows. Choose a path in
$\tilde K$ joining two such leaves. Attach a 2-cell $[0,1]\times
[0,1]$ by gluing $[0,1]\times 0$ to the path. Map the other 3 boundary
components to the arc $[x,y]\subset T$ (point if $x=y$). Then 
extend $f$ to the 2-cell in the same way as when constructing a resolution.
Finish the construction by attaching the whole orbit of 2-cells and
extending to preserve equivariance. Slight care and subdivisions may
be necessary to stay in the simplicial category. In the end, we have another
resolution $f':\tilde K'\to T$ and a factorization
$f=f'\rho$, where $\rho:\tilde K\to\tilde K'$ is equivariant, sends leaves
to leaves, and the images in $\tilde K'$ of the original pair of leaves
are joined by a path whose measure is equal to the distance
between $x$ and $y$.

Continuing in this fashion, we can construct resolutions that more
and more faithfully reflect the nature of $T$.

\begin{prop} \label{resolution}
Assume that a finitely presented group $G$ is acting
by isometries on an \rt\ $T$ and the action is non-trivial and minimal.
For any finite collection $Y\subset T$ of points in $T$ and any finite
collection $G_0\subset G$ of group elements there is a resolution
$f:\tilde K\to T$ ($K$ depends on $Y$ and $G_0$) and a collection of
points $Y'\subset \tilde K$ such that $\tilde f$ induces a bijection
between $Y'$ and $Y$ and for any $a,b\in Y'$ and any $\gamma,\delta\in
G_0$ there is a path $\alpha$ from $\gamma(a)$ to $\delta(b)$ whose
measure is equal to $d(f(\gamma(a)),f(\delta(b)))$.\qed\end{prop}

For example, if $x\in T$ and $H$ is a finitely generated subgroup of
the stabilizer $Stab(x)$, then we can construct a resolution $f:\tilde
K\to T$ such that $f$ sends a leaf or a complementary component $D$ to $x$
and $h(D)=D$ for all $h\in H$.

On the other hand, using the construction outlined in Example
\ref{tracks}, one can show that there are examples of free actions of
the free group $F_3$ such that every resolution is simplicial
\cite{bf:outerlimits}.

\subsection{Band complexes}
It is more
convenient to work with a special class of 2-complexes equipped with 
with measured laminations, called {\it band complexes}. 

\begin{definition} A {\it band} is the square $[0,1]\times [0,1]$
equipped with a measured lamination $C\times [0,1]$ with
measure of full support for a
compact totally disconnected set $C\subset (0,1)$.

A {\it multiinterval} $\Gamma$ is the disjoint union of closed intervals
equipped with a measured lamination $\Lambda(\Gamma)$ 
disjoint from the endpoints.

A {\it union of bands} is the space $Y$ obtained from a multiinterval
$\Gamma$ by attaching a collection of bands. Each band $[0,1]\times
[0,1]$ is attached via an embedding $\phi:[0,1]\times \{ 0,1\}\to \Gamma$
such that $\phi^{-1}(\Lambda(\Gamma))=C\times\{0,1\}$ and such that
$\phi$ is measure-preserving. The measured lamination $\Lambda(\Gamma)$
pieces together with the measured laminations on the bands to produce
a measured lamination $\Lambda(Y)$ on $Y$.

A {\it band complex} is the space $X$ obtained from a union of bands
$Y$ by successively attaching 0-, 1-, and 2-cells (with PL attaching
maps) so that
\begin{itemize}
\item There is a neighborhood of $\Lambda(Y)$ disjoint from the images of
all attaching maps.
\item The images of attaching maps of 1-cells are contained in $\Gamma$.
\end{itemize}
The band complex $X$ is equipped with the induced measured lamination
$\Lambda=\Lambda(X)$.
\end{definition}

\begin{example} \label{conversion} 
Let $X$ be the hyperbolic surface of Example \ref{hypsurf}.
Each triangle in $X$ intersects the lamination $\Lambda$ in a
collection of geodesic arcs, each spanning between two sides. Thus
these arcs fall into at most 3 families according to which two sides
they intersect. We can view $X$ as a band complex as follows. The
multi-interval $\Gamma$ is obtained from the 1-skeleton by removing
small disks around each vertex. Each triangle gives rise to at most 3
bands, one for each family of geodesic arcs. The vertices are the
0-cells, there are two 1-cells for each edge of the triangulation,
connecting an endpoint to $\Gamma$. Finally, a triangle of the most
interesting type (intersecting $\Lambda$ in 3 families of arcs) gives
rise to four 2-cells, three corner triangles, and a central
hexagon. Simpler triangles give rise to fewer 2-cells.
\end{example}

\begin{definition} Let $X$ be a band complex and assume that $\pi_1(X)$
is acting on an \rt\ $T$. An equivariant map $f:\tilde X\to T$ is a {\it
resolution} (or an {\it exact resolution}) if there is a triangulation of $X$
so that $f$ is a resolution (or an exact resolution) in the sense of
section \ref{resolutions}.\end{definition}

\section{Rips machine}

\subsection{Moves on band complexes}

Building on the work of Makanin \cite{gm:equations} and Razborov
\cite{ar:equations}, Rips has devised a ``machine'' that transforms any
band complex into a ``normal form''. The reference for this section
is \cite{bf:stable}. Here we only outline some aspects of the Rips
machine.

There is a list of 6 moves M0-M5 that can be applied to a band complex. The
complete list is in section 6 of \cite{bf:stable}. These moves are
analogs of the elementary moves in simple homotopy theory, but
they respect the underlying measured lamination. If a band complex
$X'$ is obtained from a band complex $X$ by a sequence of these moves,
then the following holds.

\begin{itemize}
\item There are maps $\phi:X\to X'$ and $\psi:X'\to X$ that induce an
isomorphism between fundamental groups and preserve measure.
\item If $f:\tilde X\to T$ is a resolution, then the composition
$f\tilde\psi:\tilde X'\to T$ is also a resolution, and if $g:\tilde X'\to T$
is a resolution, then so is $g\tilde\phi:X\to T$. 
\item $\phi$ and $\psi$ induce a 1-1 correspondence between the
minimal components of the laminations on $X$ and $X'$.
\item $\tilde \phi$ and $\tilde \psi$ induce quasi-isometries between the
leaves of the laminations in $\tilde X$ and $\tilde X'$.
\end{itemize}

By way of illustration, we describe one of the moves, namely (M5). An
arc $J\subset \Gamma$ is said to be {\it free} if the endpoints of $J$
are in the complement of $\Lambda$, $J$ has positive measure, and it
intersects only one attaching region of a band. A free subarc $J$ is
said to be a {\it maximal} free subarc if whenever $J'\supset J$ is a
free subarc, then $J'\cap\Lambda=J\cap\Lambda$. 

Assume that $J$ is a maximal free subarc and that $J$ is contained in
the attaching region $[0,1]\times 0$ of a band $B=[0,1]\times
[0,1]$. The move (M5) consists of collapsing $J\times [0,1]$ to
$J\times 1\cup Fr\ J\times [0,1]$. Typically, the band $B$ will be
replaced by two new bands, but if $J$ contains one or both endpoints
of the attaching region $[0,1]$, then $B$ is replaced by 1 or 0 bands.
Attaching maps of relative 1- and 2-cells whose images intersect $int\
J\times [0,1)$, can be naturally homotoped upwards.

\subsection{The
classification theorem}

For simplicity, we will assume that $X$ is a band complex and
$f:\tilde X\to T$ is a resolution of an \rt\ $T$ on which $G=\pi_1(X)$
is acting, and that the action has trivial arc stabilizers. This
assumption is not necessary, but it dramatically simplifies the
statements. The reason is that in this case it is always possible to
remove annuli from a band complex.  Imagine a band complex $X$ that
contains as a subcomplex an annulus $[0,1]\times S^1$ which is thought
of as a single band with top and bottom attached to the same arc. If
the measure of the arc $[0,1]\times p$ is positive, then the element
of the fundamental group corresponding to the loop $q\times S^1$ fixes
an arc in $T$ and is therefore trivial.  We can then collapse the
annulus to the arc and replace $X$ by the resulting complex $X'$ (this
is move (M1)).

We will also assume that $\pi_1(X)$ is torsion-free.

\begin{thm}[Rips,\cite{bf:stable}] \label{rips}
Let $X$ be a band complex such that $\tilde X$ resolves an action of
the torsion-free group $\pi_1(X)$ on an \rt\ with trivial arc
stabilizers. Then $X$ can be transformed, using moves (M0-M5), to
another band complex $X'$ with the following properties. For each
minimal component $\Lambda_i'$ of the lamination $\Lambda'$ on $X'$
there is a subcomplex $X_i'$ of $X'$ that intersects $\Lambda'$ in
$\Lambda_i'$ and these subcomplexes are pairwise disjoint. All
inclusions $X_i'\hookrightarrow X'$ are $\pi_1$-injective, and so are
all inclusions from a component of the frontier $Fr(X_i')$ into
$X_i'$. Each $X_i'$ is of one of the following 3 types:
\begin{itemize}
\item {\it Surface type:} $X_i'$ is a compact surface with negative
Euler characteristic and $\Lambda_i'$ is a (filling) geodesic measured
lamination (with respect to a hyperbolic structure on the surface).
Each component of $Fr(X_i')$ is either a point or a boundary component
of the surface.
\item {\it Toral type:} $X_i'$ is the 2-skeleton of the torus from
Example \ref{torus} with the induced lamination. Each component of
$Fr(X_i')$ is a point.
\item {\it Thin type:} This type does not have a standard model. Its main
feature is that it can be arranged that $X_i'$ contains an arbitrarily thin
band (i.e. with attaching regions of small measure) that intersects
the rest of $X'$ only in the two attaching regions. See more on this below.
\end{itemize}
\end{thm}

If $\pi_1(X)$ is not torsion-free, the theorem still holds
provided that in the surface and toral types we allow for a finite number of
cone-type orbifold points.

Traditionally, in terms of the dual tree (or the associated
pseudogroup), surface type is called ``interval exchange'', and
toral type is called ``axial''. Similarly, thin type is also
called ``Levitt type'', in honor of G. Levitt \cite{gl:thin} who
discovered and extensively studied this kind of a pseudogroup. Thin
type has also been called ``exotic''. We chose names that reflect
the nature of the band complex, not the dual tree or the pseudogroup.

\subsection{Thin type band complexes}

We now describe band complexes of thin type in more detail.
Suppose $X$ is a band complex whose lamination $\Lambda$ is minimal.
If there are maximal free arcs, choose one and perform the collapsing move
(M5) described above to obtain another band complex $X_1$. If $X_1$
contains maximal free arcs, choose one and collapse, etc. This
process stops when there are no more free arcs. Note that a collapse might
produce new free arcs. This process is called Process I in \cite{bf:stable}.

\begin{definition} The band complex $X$ is of {\it thin type} if it
is equivalent under moves (M0-M5) to a complex for which the
collapsing procedure never ends.\end{definition}

There is a concrete example of a thin band complexes in section 10 of
\cite{bf:stable}.  This example has the additional feature that after
each collapse the resulting band complex is a scaled down version of
the original. R. Martin \cite{reiner:thin} has studied the relation between
the ``periodicity'' of the sequence of collapses and unique ergodicity
of the underlying lamination.
The band complexes associated to the ``interesting'' pseudogroups in
\cite{gl:thin} are thin.

If we focus on a single band in $X$, then under the collapsing process
this band will get subdivided into more and more bands with
arbitrarily small transverse measure (so these bands are thin, thus
the name). In particular, eventually there will be bands whose
interiors are disjoint from the attaching regions of the relative
2-cells. Such bands are called {\it naked bands}. A naked band induces
a free product decomposition of $\pi_1(X)$ by cutting along an arc in
the band that separates the two attaching regions. That this
decomposition is non-trivial is the content of Proposition 8.13 of
\cite{bf:stable}.

As the reader will learn from section 7, in applications one
frequently assumes that the underlying group $\pi_1(X)$ is freely
indecomposable, and then the situation simplifies considerably as
there can be no thin components in resolving band complexes (of
course, assuming the arc stabilizers of $T$ are trivial). Similarly,
when one is concerned with hyperbolic groups, there can be no toral
components.

One can also make a study of quasi-isometry types of leaves for a thin
type lamination. Generic leaves are quasi-isometric to 1-ended trees,
and in addition there are uncountably many leaves quasi-isometric to
2-ended trees. For details see Proposition 8.13 of \cite{bf:stable}
and, independently, Gaboriau \cite{gaboriau:orbits}. Of course, this
is to be contrasted with the surface and toral types where the leaves
are quasi-isometric to Euclidean space (of dimension 1 and $>1$
respectively).

\subsection{Remarks on the proof of the
classification theorem} 

We now briefly describe the proof of Theorem \ref{rips}. For details,
see section 7 of \cite{bf:stable}. As mentioned in the introduction,
there is an alternative approach developed by Gaboriau, Levitt, and
Paulin.

There is an algorithm for transforming a given band complex (say with
a minimal lamination) to another one. When there is a free arc, we
collapse from a maximal free subarc as described above. This is
Process I. When there are no free subarcs, one performs a {\it
sliding} move, called Process II. These are to be repeated producing a
sequence of band complexes. There is a notion of complexity
(non-negative half-integer valued). The moves never increase
complexity, and whenever Process II is followed by Process I
(i.e. whenever free subarcs disappear after a collapse) the complexity
strictly decreases. It follows that eventually only Process I or only
Process II is performed. In the first case the band complex is of thin
type, and in the second one argues that it is of surface or of toral
type.

\section{Stable actions on \rt s}

Here is a sample statement that illustrates how the Rips machine can be
applied to understand the structure of a finitely presented group that
is acting on an \rt.

\begin{thm} \label{easy}
Suppose that a torsion-free finitely presented group $G$ is acting
non-trivially on an \rt\ $T$ by isometries and that all arc
stabilizers are trivial. Then one of the following holds.
\begin{itemize}
\item $G$ splits as a non-trivial free product. In this case one can study
the free factors by examining the induced action on $T$. Either this
theorem can be applied to a given factor or this factor is a point stabilizer
in $T$.
\item $G$ is a free abelian group.
\item $G$ is the fundamental group of a 2-complex $X$ that contains as
a subcomplex a compact connected surface $S$ of negative Euler
characteristic and $S\cap \overline{X\setminus S}$ is contained in
$\partial S$. Inclusion induced homomorphism $\pi_1(S)\to\pi_1(X)=G$
is injective and each boundary component of $S$ corresponds to an
elliptic isometry in $T$. There is a filling geodesic measured
lamination with measure of full support on $S$ disjoint from $\partial
S$ (``filling'' means that each complementary component is either
simply-connected or it contains a boundary component and its
fundamental group is $\Z$).
\end{itemize}
\end{thm}

\begin{proof}
We may assume that the action is minimal by passing to the minimal
subtree, as in Proposition \ref{countably many axes}. Let $K$ be a
finite complex with $\pi_1(K)=G$ and choose a resolution $f:\tilde
K\to T$ as in section 4.3. Convert $K$ to a band complex $X$ in the same
manner as in Example \ref{conversion}. Then apply Theorem \ref{rips}
to replace $X$ by a band complex $X'$ in ``normal form''.

If some component $X_i'$ is of thin type, then $X_i'$ can be assumed to
contain a naked band and hence $G=\pi_1(X')$ splits as a non-trivial free
product, so the first possibility holds.

If some component $X_i'$ is of toral type, then either we obtain
a free product decomposition of $G$ using one of the points in $Fr(X_i')$
(so the first possibility holds),
or $G=\pi_1(X_i')$ is free abelian, so the second possibility holds.

If some component $X_i'$ is of surface type, then the third possibility
holds.

Finally, if the lamination on $X'$ is simplicial, then $G$ acts on the
simplicial tree dual to this lamination. The action can only be
``freer'' than the original action, so from Bass-Serre theory we
conclude that the first possibility holds.
\end{proof}

E. Rips presented a proof of the following theorem at the conference
at the Isle of Thorns in the summer of 1991. It answers affirmatively
the conjecture of Morgan and Shalen.

\begin{thm}[Rips] \label{free}
If $G$ is a finitely presented group that acts freely by isometries
on an \rt, then $G$ is the free product of free abelian groups and
closed surface groups.
\end{thm}

\begin{proof}
Decompose $G$ into the free product of a free group and freely
indecomposable factors, and apply Theorem \ref{easy} to each freely
indecomposable factor.
\end{proof}

As indicated earlier, the methods generalize to stable actions.
The following is stated as Theorem 9.5 in \cite{bf:stable}.

\begin{thm} \label{useful}
Let $G$ be a finitely presented group with a stable action on
an \rt\ $T$. Then either
\begin{itemize}
\item $G$ splits over an extension $E$-by-cyclic where $E$ fixes an arc
of $T$, or
\item $T$ is a line and $G$ splits over an extension of the kernel of
the action by a free abelian group.
\end{itemize}
\end{thm}

The structure of the group $G$ obtained in theorems above very much
depends on the choice of the resolution. Imagine taking a sequence
of finer and finer resolutions of the given stable action, as in the
discussion preceding Proposition \ref{resolution}. Each band complex
in the sequence gives rise to a splitting of $G$ (more precisely, to
a graph of groups decomposition of $G$).

\vskip .2cm
\noindent {\it Question.} Does the sequence of splittings of $G$
``stabilize''? In other words, is there a structure theorem for $G$
that does not depend on the choice of a resolution, but only on the
tree?

\vskip .2cm
\noindent {\it Question.} Does Theorem \ref{useful} hold if ``finitely
presented'' is replaced by ``finitely generated'' in the hypotheses?
\vskip .2cm

Theorem \ref{free} holds in the setting of finitely generated groups.
Zlil Sela answered the above two questions affirmatively in the case that
the action has the additional property that the stabilizers of tripods
are trivial. This important case often arises in applications.

The structure of the
group acting on an \rt\ without the assumption of stability is still
very much a mystery.

\vskip .2cm
\noindent {\it Question.} If a finitely presented group $G$ admits a
non-trivial isometric action on an \rt, does it also admit a
non-trivial action on a simplicial tree (i.e. does it admit a
non-trivial splitting)?
\vskip .2cm

The answer is affirmative if $G$ is a 3-manifold group by
the work of Morgan and Shalen (Proposition 2.1 of \cite{ms:valuations3}).

\section{Applications}

I will now outline a number of applications of the theory of
$\R$-trees.  The technique can naturally be used in proofs of
finiteness and compactness theorems. Surprisingly, as shown by the
work of Zlil Sela, $\R$-trees can also be used to derive various
structure theorems in group theory. It is impossible to cover all
applications to date, so it seems reasonable to restrict this
exposition to outlines of the most typical and the most striking
applications. There is no discussion of the work of Rips and Sela on
JSJ decompositions of finitely presented groups (see \cite{rs:jsj}),
in part because in the meantime simpler proofs of more general
theorems have been found \cite{ds:jsj},\cite{fp:jsj}. There is no
doubt, however, that the intuition coming from the theory of \rt s played
a key role in this discovery.

\subsection{Compactifying spaces of geometric structures}

Topologists became interested in $\R$-trees with the work of Morgan
and Shalen \cite{ms:valuations1} that shed new light and
generalized parts of Thurston's Geometrization Theorem. If $M$ is a
closed oriented $n$-manifold, then having a hyperbolic structure on $M$ is
equivalent to having a discrete and faithful representation
$\pi_1(M)\to Isom_+\H^n$ into the orientation-preserving isometry
group of the hyperbolic $n$-space, up to conjugation in $Isom_+\H^n$.
For $n=2$ and $M$ of genus $g\geq 2$ the space
$$Hom_{DF}(\pi_1(M),Isom_+\H^2)/conj$$ of hyperbolic structures on $M$
is the {\it Teichm\" uller space} of $M$. It is known that this space
is homeomorphic to Euclidean space of dimension $6g-6$. The
automorphism group $Aut(M)$ of $M$ (homeomorphisms of $M$ modulo
isotopy, also known as the {\it mapping class group} of $M$) naturally
acts on it (with finite isotropy groups), so the Teichm\" uller space
is useful in the study of $Aut(M)$ as it plays the role of the
classifying space.

An important ingredient of Thurston's theory of surface automorphisms
\cite{wt:diffeos},\cite{flp:surfaces} is his construction of an
equivariant compactification of the Teichm\" uller space.  An ideal
point is represented by a transversely measured geodesic lamination on
$M$ (measures that differ by a multiple are equivalent).

From the point of view of $\R$-trees, the construction of this
compactification comes from the Compactness Theorem (see section 3.5). 
An ideal point is
represented by a non-trivial and minimal isometric action of $\pi_1(M)$ on
an $\R$-tree, with homothetic actions
considered equivalent. Further, from Proposition \ref{small} we see
that the arc stabilizers are cyclic. Recall that an action of
$\pi_1(M)$ on an $\R$-tree is small if it is minimal, does not have
global fixed points, and all arc stabilizers are cyclic.

That the two approaches are equivalent follows from the following
result of Skora \cite{rs:splittings}.

\begin{thm}[Skora \cite{rs:splittings}] \label{skora}
If $M$ is a closed hyperbolic surface, then any small
action of $\pi_1(M)$ on an $\R$-tree is dual to a unique measured
geodesic lamination on $M$.\end{thm}

\begin{proof} The proof using the Rips machine is considerably simpler
than the original proof. We focus on the special case when the action
is free; the general case is similar. Let $f:\tilde X\to T$ be
a resolution of the action. Since $\pi_1(M)$ does not contain $\mathbb
Z\times\mathbb Z$, $X$ cannot have toral components, and since it is
not freely decomposable, $X$ cannot have simplicial or thin
components, and in fact $X$ must have a single surface component (see
Theorem \ref{easy}). Thus $X$ can be taken to be a closed surface
equipped with a measured geodesic lamination that fills the
surface. To finish the proof, we need to argue that $f$ is an exact
resolution. If not, $f$ factors as $f=gh$ through another resolution
$g:\tilde X'\to T$, where $\tilde h:\tilde X\to \tilde X'$ is
equivariant and sends leaves to leaves.  By the same argument as
above, $X'$ can be taken to be a closed surface with a filling
measured geodesic lamination. Thus $h:X\to X'$ is a homotopy
equivalence that sends leaves to leaves and locally preserves the
transverse measure. In the universal cover (identified with the
hyperbolic plane), distinct leaves diverge from each other (in at
least one direction), and therefore $\tilde h$ cannot send distinct
leaves to the same leaf. Since $\tilde h$ induces a homeomorphism
between the circles at infinity and a lamination is determined by
the pairs of endpoints at infinity of its leaves, it follows that $h$
can be taken to be a homeomorphism, showing that $f$ is an exact
resolution.
\end{proof}

Theorem \ref{skora} plays a prominent role in J.-P. Otal's proof
\cite{otal:hyperbolization} of Thurston's Double Limit Theorem, which
in turn is a key ingredient in the proof of the Hyperbolization
Theorem for 3-manifolds that fiber over the circle.

In dimensions $n>2$ the celebrated Rigidity Theorem of Mostow states
that the space of hyperbolic structures on a closed manifold $M^n$ has
at most one point, and the construction using the Compactness Theorem
is not particularly exciting in that case. However, it is important in
Thurston's proof of the Geometrization Theorem to study the space
$$Hom_{DF}(G,Isom_+\H^n)/conj$$ where $G$ is the fundamental group of
a compact 3-manifold (with boundary) and $n=3$. In particular,
Thurston needed the fact that this space is compact when the
3-manifold is irreducible, aspherical, acylindrical, and atoroidal. In
group-theoretic terms, this means that $G$ is torsion-free and does
not split over $1,\Z$, or $\Z^2$. 

\begin{thm} \cite{bf:stable}
Suppose $G$ is finitely presented, not virtually abelian, 
and does not split over a
virtually abelian subgroup. Then the space
$$Hom_{DF}(G,Isom_+\H^n)/conj$$
of homotopy hyperbolic structures on $G$ is compact.\end{thm}

\begin{proof} If the space is not compact, there is a sequence going to
infinity. The Compactness Theorem provides a small action of $G$ on an
$\R$-tree. Theorem \ref{useful} then implies that $G$ splits over a
virtually abelian subgroup. (Recall that a discrete group of isometries
of $\H^n$ is either virtually abelian or it contains $F_2$.)\end{proof}

This theorem generalizes earlier
work of Thurston, Morgan-Shalen, and Morgan.

\subsection{Automorphism groups of word-hyperbolic groups}

It is the fundamental observation of F. Paulin \cite{fp:outer} that
$\R$-trees arise also in the coarse setting of word-hyperbolic
groups in the presence of infinitely many automorphisms of the group. 
The second part of the proof of the following theorem follows
from the Rips machine.

\begin{thm} \label{paulin}
Suppose $G$ is a word-hyperbolic group such that $Out(G)$ is
infinite. Then $G$ splits over a virtually cyclic subgroup.\end{thm}

\begin{proof} Let $f_j:G\to G$ be an infinite sequence of pairwise
non-conjugate automorphisms. Each $f_j$ produces an isometric action
$\rho_j$ of $G$ on its Cayley graph by sending $g\in G$ to the left
translation by $f_j(g)$. The Compactness Theorem
provides an action of $G$ on an $\R$-tree $T$. The arc stabilizers of
this action are small by Proposition \ref{small}, so the claim follows
from Theorem \ref{useful}.\end{proof}

We will now assume that $G$ is a torsion-free word-hyperbolic group.
It is an open question whether every word-hyperbolic group has a
torsion-free subgroup of finite index (or even whether it is
residually finite). It is known that there are only finitely many conjugacy
classes of finite order elements \cite{mg:hyperbolic}.

For torsion-free $G$ the converse of Theorem \ref{paulin} holds. If
$G$ splits as a free product $G=A*B$ with $A$ and $B$ nontrivial
(infinite!), and if one, say $A$, is nonabelian, then for a fixed
nontrivial $a\in A$ the automorphism $f:G\to G$ that restricts to
identity on $B$ and to conjugation by $A$ represents an element of
infinite order in $Out(G)$.  The remaining case is $G=\Z*\Z=F_2$, but
$Out(F_2)=GL_2(\Z)$ has many elements of infinite order. If $G$ splits
over $\Z$, say as $G=A*_{C}B$, with $A\neq C\neq B$ and $C=<c>$
infinite cyclic, then there is a {\it Dehn twist} automorphism $f:G\to
G$ that restricts to identity on $B$ and to conjugation by $c$ on
$A$. This represents an element of infinite order as long as $A$ is
nonabelian. If $A$ is abelian, then reverse the roles of $A$ and
$B$.  (Not both $A$ and $B$ could be abelian, since then $G$ would
not be word-hyperbolic.)  Finally, if $G$ splits as $G=A*_C$ with
$C=<c>$ infinite cyclic, then the automorphism ({\it Dehn twist}) that
restricts to the identity on $A$ and sends the ``stable letter'' $t$
to $tc$ has infinite order in $Out(G)$.

The above paragraph suggests the following sharpening of Theorem
\ref{paulin}, at least for torsion-free, freely indecomposable
word-hyperbolic groups.  Call the subgroup of $Aut(G)$ generated by
all inner automorphisms and all Dehn twists (with respect to all
possible splittings over infinite cyclic subgroups) the {\it Internal
Automorphism Group}, denoted $Int(G)$. It is a normal subgroup of
$Aut(G)$.  Note that the celebrated theorem of Dehn \cite{dehn:twists}
(see also \cite{wbrl:twists}) that the mapping class group of a closed
orientable surface is generated by Dehn twists can be interpreted as
saying $Int(G)=Aut(G)$ where $G$ is the fundamental group of the
surface.  If the surface is allowed to be non-orientable and to have
boundary, then the subgroup of the automorphism group (i.e. the
homeomorphism group modulo isotopy) generated by Dehn twists has
finite index.

\begin{thm}[Rips-Sela \cite{rips-sela:gafa}] \label{internal}
If $G$ is a torsion-free, freely indecomposable word-hyperbolic group, then
the Internal Automorphism Group has finite index in $Aut(G)$.\end{thm}

The proof introduces a new idea, the {\it shortening argument}.

\begin{proof}
Fix a finite generating set $\{\gamma_1,\cdots,\gamma_k\}$ for $G$
which is closed under taking inverses and for $f\in Aut(G)$ define
$$d(f)=\max_{1\leq i\leq k}||f(\gamma_i)||$$ 
where $||\cdot ||$ denotes the word length. In each coset of $Int(G)$ in
$Aut(G)$ choose an automorphism $f$ with minimal $d(f)$. Assuming that
there are infinitely many cosets, we have an infinite sequence of
automorphisms $f_1,f_2,\cdots\in Aut(G)$ that represent distinct
cosets of $Int(G)$ and each minimizes the function $d$ in its
coset. As in the proof of Theorem \ref{paulin}, we view each $f_j$ as
giving an action $\rho_j$ of $G$ on its Cayley graph. Note that $1\in
G$ is centrally located for $\rho_j$ (or else composing $f_j$ with an
inner automorphism would produce a representative of the same coset
with smaller $d$).

After passing to a subsequence, we obtain a limiting action $\rho$ of $G$
on an $\R$-tree $T$. We will examine this action and argue that for
large $j$ the automorphism $f_j$ can be composed with a Dehn twist in such
a way that $d$ is reduced.

Let $X$ be a finite band complex and $\phi:\tilde X\to T$ a resolution
of $\rho$. Choose a basepoint $*\in \tilde X$ that maps to the
basepoint in $T$. We can arrange (see Proposition \ref{resolution})
that the distance in $\tilde X$ between $*$ and $\gamma_i(*)$ equals
the corresponding distance in $T$.

Since $G$ is word-hyperbolic and so does not contain $\Z\oplus\Z$, $X$
cannot have any toral components. Similarly, $X$ cannot have any thin
components, as we are assuming that $G$ is freely indecomposable. Therefore,
all components of $X$ are of simplicial or surface type, and the simplicial
pieces have infinite cyclic edge stabilizers.

Let us consider the two extreme cases. First suppose that $X$ is a
closed surface with a filling geodesic-like measured lamination. It is
a fact of surface theory that there is a homeomorphism $h:X\to X$
fixing the basepoint, which can be taken to be a product of Dehn
twists, such that the measure of each $h([\gamma_i])$ is arbitrarily
small.  This fact can be proved by ``unzipping'' the band complex
(i.e. the ``train-track'', see \cite{flp:surfaces}) until the bands
are arbitrarily thin and taking for $h$ a homeomorphism that sends
thick bands to thin bands. Now $f_j\pi_1(h)$ is a ``shorter''
representative of the coset $f_jInt(G)$ for sufficiently large $j$, 
a contradiction.

Now suppose that $X$ is simplicial. Let $T'$ be the simplicial tree
dual to $\tilde X$. For notational simplicity we assume that $T'/G$ is
a single edge, corresponding to an amalgamated product decomposition
of $G$ (over $\mathbb Z$).  The basepoint in $\tilde X$ corresponds to a
vertex $v$ in $T'$. Consider an edge $e$ of $T'$ that has $v$ as an
endpoint. The stabilizer of $e$ is infinite cyclic. Say $c\in G$
generates this stabilizer. Let $A$ denote the stabilizer of $v$ and
$B$ the stabilizer of the other endpoint of $e$, so that
$G=A*_{<c>}B$. Also, without loss of generality we can assume that
the length of $e$ is 1. The distance between $v$ and $\gamma_i(v)$ in $T'$
(equivalently, in $T$) is the minimal $m_i$ such that $\gamma_i$ is
the product of the form $a_0b_1a_2\cdots $ of total length $m_i+1$ with the 
$a$'s in $A$ and the $b$'s in $B$.


Fix a large $j$
and consider the translates of the basepoint $1\in G$ under the
generators $\gamma_i$ with respect to the representation $\rho_j$.
After rescaling by the constant $d_j=d(f_j)$, the word-metric on $G$ restricted
to this finite set is close to the metric induced from $T$ (or $T'$) by
restricting to the translates of $v$ by the generators.


Let $b\in B$ be one of the $b$'s occurring in the above representations
of the $\gamma_i$'s. The axis of $f_j(c)$ in $G$ and the geodesic
joining 1 and $f_j(b)$ are within $10\delta$ for a length of
about $d_j$ and the translation length of $f_j(c)$ is $<<d_j$.
Replace $c$ by $c^{-1}$ if necessary so that $f_j(c)$ translates from
$f_j(b)$ towards 1. Choose the (positive) power $m$ so that
$f_j(c)^m$ translates $f_j(b)$ about halfway towards 1. 

We now claim that precomposing $f_j$ by the $m^{th}$ power $h$ of the
Dehn twist that fixes $A$ and conjugates $B$ by $c$ has the effect of
shortening the representative of the coset $f_jInt(G)$. Indeed, write
$\gamma_i=a_0b_1a_2\cdots$ so that
$f_jh(\gamma_i)=f_j(a_0)f_j(b_1)^{f_j(c)^m}f_j(a_2)\cdots$. The
distance between $1$ and $f_jh(\gamma_i)$ can be estimated in the usual way:
\begin{eqnarray*}
d(1,f_jh(\gamma_i))\leq &
d(1,f_j(a_0))+d(f_j(a_0),f_j(a_0)f_j(b_1)^{f_j(c)^m})+&\\
&d(f_j(a_0)f_j(b_1)^{f_j(c)^m},f_j(a_0)f_j(b_1)^{f_j(c)^m}f_j(a_2))+\cdots=&\\
&d(1,f_j(a_0))+d(1,f_j(b_1)^{f_j(c)^m})+d(1,f_j(a_2))+\cdots&\end{eqnarray*} 
The terms
of the form $d(1,f_j(a))$ are small compared to $d_j$ (the ratio goes
to 0), and the terms $d(1,f_j(b_k)^{f_j(c)^m})$ are approximately
$d_j/2$. Thus the distance $d(1,f_jh(\gamma_i))$ is estimated above by
about $m_id_j/2$, and this is much less than $d(1,f_j(\gamma_i))$
(which is about $m_id_j$).

The general case (when $T'$ has perhaps more than one orbit of edges, or
when $X$ has both surface and simplicial components) is dealt with
in the same way; only notation is more involved.
\end{proof}

A version of the theorem can be proved for torsion-free word-hyperbolic
groups that are free products using the classical theory of automorphisms
of free products \cite{fr:free1},\cite{fr:free2}.

The same method has other applications. Recall that a group $G$ is
{\it co-Hopfian} if every injective endomorphism $G\to G$ is
surjective. Nontrivial free products are never co-Hopfian. For our
purposes, the group $\Z$ is not freely indecomposable (it splits over
the trivial group).

\begin{thm}[Sela \cite{sela:cohopf}] 
Every freely indecomposable word-hyperbolic group is
co-Hopfian.\end{thm}

\begin{proof} Let $Inj(G)$ denote the semi-group
of injective endomorphisms of $G$. The idea is to follow the above
argument and show that $Aut(G)$ has finite index in $Inj(G)$. The only
difference with the situation $Int(G)\subset Aut(G)$ is that $Inj(G)$
is not a group and $Aut(G)$ is not normal in $Inj(G)$, but those
features were never used.  Finally, note that if $Aut(G)$ has finite
index in $Inj(G)$, then a nontrivial power of every $f\in Inj(G)$ is
an automorphism, and so $Inj(G)=Aut(G)$.
\end{proof}

Recall that a group $G$ is {\it Hopfian} if every surjective
endomorphism $G\to G$ is an isomorphism. Z. Sela has announced the following
result \cite{sela:hopf}:

\begin{thm} Every torsion-free word-hyperbolic group is Hopfian.
\end{thm}

The proof uses more elaborate ideas and will not be outlined here.

\begin{thm}[Gromov \cite{mg:hyperbolic}, Sela \cite{sela:cohopf}] Let
$\Gamma$ be a finitely presented torsion-free freely indecomposable
non-cyclic group and $G$ a word-hyperbolic group. Then there are only
finitely many conjugacy classes of subgroups of $G$ isomorphic to
$\Gamma$.\end{thm}

\begin{proof} First consider the simple case when $\Gamma$ does not admit
any splittings over $\mathbb Z$. Then we argue that there can be only
finitely many conjugacy classes of monomorphisms $f:\Gamma\to G$. For
suppose that there are infinitely many. Let $f_1,f_2,\cdots:\Gamma\to
G$ be an infinite sequence of pairwise non-conjugate monomorphisms. We
thus get a sequence of actions $\rho_i$ of $\Gamma$ on $G$:
$\rho_i(\gamma)$ acts by left translation by $f_i(\gamma)$.  By
conjugating each $f_i$ we may assume that $1\in G$ is centrally
located with respect to each $\rho_i$ (and with respect to a fixed
finite generating set for $\Gamma$). Now pass to a subsequence and
obtain an action of $\Gamma$ on an \rt. As before, this action induces
a splitting of $\Gamma$ over $\Z$.

\hyphenation{mono-morphism}

If $\Gamma$ admits a splitting over $\Z$, then we could precompose a
given monomorphism $\Gamma\to G$ by automorphisms of $\Gamma$ and
obtain an infinite sequence of non-conjugate monomorphisms $\Gamma\to
G$. This phenomenon is precisely what the shortening argument is
designed to handle. Given a monomorphism $\Gamma\to G$, conjugate it
by an element of $G$ and precompose by an automorphism of $\Gamma$ so
as to make $1\in G$ centrally located and to make the maximal
displacement of 1 smallest possible. Now the claim is that there can
be only finitely many such minimizing monomorphisms. The proof of the
claim is analogous to the proof of Theorem \ref{internal}. If there
are infinitely many such monomorphisms, consider the limiting tree and
use it to construct an automorphism $h:\Gamma\to\Gamma$ that can be
used to shorten representations $\rho_j$ for large $j$.
\end{proof}

\subsection{Fixed subgroup of a free group automorphism}

Let $F_n$ be the free group of rank $n$ and $f:F_n\to F_n$ an automorphism.
Recall that $F_n$ contains free subgroups of infinite rank. The following
theorem was conjectured by Peter Scott.

\begin{thm}[Bestvina-Handel] \label{scott}
The rank of the subgroup $Fix(f)$ of elements of $F_n$ fixed by $f$
is at most $n$.\end{thm}

The proof in \cite{bh:tracks} does not use the theory of $\R$-trees. Z. Sela
\cite{sela:scott} and D. Gaboriau-G. Levitt-M. Lustig \cite{gll:scott}
found a simpler argument using $\R$-trees. We now outline their ideas.

First, for $k=1,2,\cdots$ let $g_k$ be an automorphism conjugate to
$f^k$ such that $1\in F_n$ is centrally located with respect to the
representation $\rho_k$ that to $\gamma\in F_n$ associates left translation
$F_n\to F_n$ by $g_k(\gamma)$. This conjugation is necessary in order
to apply the Compactness Theorem, but of course $Fix(g_k)$ is in
general different from $Fix(f)$. It is therefore more natural to
consider elements of $F_n$ fixed up to conjugacy. If $f$ has finite
order as an outer automorphism, the rescaling constants remain
bounded. Such automorphisms were handled by Culler \cite{culler:finite} who
showed that the fixed subgroup is either cyclic or a free factor. For
non-periodic automorphisms, we analyze the action of $F_n$ on the
$\R$-tree $T$ obtained as the limit of a subsequence of
representations $\rho_k$ above.

The key observation is that any $\gamma\in F_n$ which is fixed up to
conjugacy by $f$ is elliptic in $T$. Indeed, the translation length of
$\gamma$ can be computed as the limit of ratios
$$\frac{{\rm translation\ length\ of}\ g_k(\gamma)}{{\rm rescaling\
factor\ for}\ \rho_k}$$ and this converges to 0 since the denominators
go to infinity, while the numerators are constant (and equal to the
length of the conjugacy class of $\gamma$). The same argument shows
that periodic conjugacy classes are elliptic in $T$ (and also those
that grow slower than the fastest growing conjugacy classes).

Second, we construct a bilipschitz homeomorphism $H:T\to T$ which is
equivariant with respect to $f$, i.e. $h(\gamma(x))=f(\gamma)(h(x))$.
This construction is due to Sela who used it extensively. He calls it
the ``basic commutative diagram''. First form the group
$G=F_n\rtimes_f\Z=<F_n,t|tgt^{-1}=f(g)>$, the mapping torus of $f$. Each
action
$\rho_k$ extends to an action
$\tilde \rho_k$ of $G$ on $F_n$ by sending $t$ to the conjugate of $f$
by the same element used to conjugate $f^k$. Of course, the extended
action is not isometric, only bilipschitz. Pass to a subsequence as usual to
obtain a bilipschitz action of $G$ on an $\R$-tree. Restricting to $F_n$
gives the discussion of the first paragraph, while $t\in G$ provides the
desired bilipschitz homeomorphism $H:T\to T$.

Third, we promote $H$ to a homothety. This is not absolutely necessary
here, but in other applications it comes handy. The following
construction is due to Paulin \cite{pa:schauder}. The Compactness
Theorem implies that the space ${\cal PED}_0$ of projectivized
nontrivial 0-hyperbolic equivariant distance functions on $F_n$ is
compact. The preimage of the closed subset ${\cal {PED}}_0^T$ of ${\cal PED}_0$
consisting of those projective classes of distance functions $d$ with
the property that
$$(x\cdot y)_T\geq (x\cdot z)_T\Rightarrow (x\cdot y)_d\geq (x\cdot
z)_d$$ in $\cal ED$ is a convex cone: If $d_1$ and $d_2$ are two
0-hyperbolic distance functions satisfying the above condition, then
$sd_1+(1-s)d_2$ is also such a distance function for $0\leq s\leq
1$. Subscripts $T$ and $d$ above indicate the metric with respect to
which $(\cdot,\cdot)$ is taken. It easily follows that ${\cal
{PED}}_0^T$ is a compact absolute retract, and therefore has the fixed
point property.  By pulling back, $H$ induces a homeomorphism of
${\cal PED}_0^T$. A fixed point of $H$ determines a new 0-hyperbolic
distance function on $F_n$ with respect to which $H$ is a
homothety. By Connecting the Dots (Lemma \ref{connecting-the-dots}),
we obtain a new tree that we continue to denote by $T$. We remark that
the new tree may not be homeomorphic to the old, but is rather
obtained from the old by collapsing some subtrees. What is important
is that arc stabilizers in the new tree are contained in the arc
stabilizers of the old tree.

Alternatively, steps 1-3 could have been avoided by quoting 
some of the theory developed in \cite{bh:tracks}. See \cite{lustig:trees},
where this alternative construction is carried out in detail.

We now arrive at the heart of the argument.

\begin{prop} \label{gl}
Assume that $F_n$ acts on an \rt\ $T$ and the action is small. Then all
vertex stabilizers of $T$ have rank $\leq n$. Further, if there is a vertex
stabilizer $V$ of rank $n$, then the action is simplicial, all edge
stabilizers are infinite cyclic, and every vertex stabilizer that is
not infinite cyclic is conjugate to $V$. \end{prop}

Before giving the proof of Proposition \ref{gl} we finish the proof of
Theorem \ref{scott}. We have seen above that each $\gamma\in Fix(f)$
is elliptic in $T$. It is an exercise to show that there is a point
$v\in T$ fixed by each $\gamma\in Fix(f)$. The ingredients are 1) the
product of two elliptic isometries of $T$ with disjoint fixed point
sets is hyperbolic, and 2) arc stabilizers of $T$ are cyclic. We may
assume that $rank(Fix(f))>1$, and then $v$ is unique. Since $f$ leaves
$Fix(f)$ invariant, equivariance forces $H$ to fix $v$. In particular,
$H$ induces an automorphism $f_v:Stab(v)\to Stab(v)$. If
$rank(Stab(v))<n$, we can apply induction on the rank and conclude
that $rank(Fix(f))=rank(Fix(f_v))<n$.  If $rank(Stab(v))=n$ and $f_v$
has finite order (as an outer automorphism) we can apply Culler's
result to conclude that $rank(Fix(f))\leq n$. If $rank(Stab(v))=n$ and
$f_v$ has infinite order, we can repeat the construction with $f_v$ in
place of $f$. We obtain a sequence of automorphisms
$f=f_0,f_v=f_1,f_2,\cdots$. We can stop when the rank of the vertex
stabilizer is $<n$ or when the restriction of the automorphism to the
vertex stabilizer has finite order. It remains to argue that the
sequence must terminate. The tree $T$ constructed above provides a
graph of groups decomposition ${\cal G}_0$ of $F_n$ with cyclic edge
groups (according to Proposition \ref{gl}). The only vertex group is
$Stab(v)$ and it has rank $n$. The next iteration provides a graph of
groups decomposition ${\cal G}_1$ of $Stab(v)$ of the same nature. We
claim that ${\cal G}_1$ can be used to refine ${\cal G}_0$, by
``blowing up'' the vertex. Indeed, this will be possible if all edge
groups of ${\cal G}_0$ are elliptic in ${\cal G}_1$. But the edge
groups of ${\cal G}_0$ are permuted (up to conjugacy) by $f$, since
the orbits of edges are permuted by $H$, and the claim follows from
the observation above that $f$-periodic conjugacy classes are elliptic
in $T$. If the sequence of automorphisms does not terminate, then
continuing in this fashion we obtain graph of groups decompositions of
$F_n$ with all edge groups cyclic, with only one vertex, and with more
and more edges. This is not possible, for example by the generalized
accessibility theorem of \cite{bf:bounding}, or better yet, by
abelianizing there can be at most $n$ edges.

\begin{proof} [Proof of Proposition \ref{gl}]
Inductively, we assume
that Proposition \ref{gl} holds for free groups of rank $<n$. If
$Stab(v)$ is contained in a proper free factor of $F_n$, then the
statement follows inductively on the rank of the underlying free
group.
\vskip .2cm
\noindent {\it Claim 1.} If $T$ is simplicial then all vertex
stabilizers have rank $\leq n$. If there is a vertex stabilizer of
rank $n$, then all other vertex stabilizers have rank 1 and all edge
stabilizers are infinite cyclic.

An edge of $T$ with trivial stabilizer induces a free factorization of
$F_n$ which implies the claim by induction. So we can assume that all
edge stabilizers are infinite cyclic. Now construct the graph of
spaces associated with the graph of groups $T/F_n$ as in
\cite{sw:topmethods}. Every vertex in $T/F_n$ is represented by a rose,
and every edge by an annulus. Since adding annuli does not change the
Euler characteristic, we see that the Euler characteristic of the
resulting space, which must be $1-n$, is equal to the sum $\sum (1-r_i)$,
where $r_1,r_2,\cdots$ denote the ranks of the vertex labels in $T/F_n$.
Since $r_i>0$ for each $i$ by assumption, Claim 1 follows.

Assume now $rank(Stab(v))<\infty$ and let $\tilde X\to T$ be a
resolution of $T$ such that a compact set $K$ in some complementary
component $D\subset X$ satisfies $im[\pi_1(K)\to\pi_1(X)]=Stab(v)$
(see Proposition \ref{resolution}). Since $F_n$ does not contain
$\Z\times \Z$ nor an extension of $\Z\times \Z$ by $\Z$, $X$ cannot
have any toral components. Likewise, if $X$ has a thin component, we
can transform it so that there is a naked band disjoint from $K$ and
we conclude that $Stab(v)$ is contained in a proper free
factor. Therefore, $X$ consists of surface and simplicial components.
If there is at least one surface component (of negative Euler
characteristic), then the Euler characteristic count of Claim 1
implies $rank(Stab(v))<n$. So assume that all components of $X$ are
simplicial, and let $T'$ be the dual (simplicial) tree. If an edge
stabilizer in $T'$ is trivial, then all vertex stabilizers in $T'$,
including $Stab(v)$, have rank $<n$. The following claim concludes the
proof in case $rank(Stab(v))<\infty$:
\vskip .2cm
\noindent {\it Claim 2.} If all edge stabilizers in $T'$ are infinite cyclic,
then $T$ is simplicial.

To prove the claim, for each primitive $a\in F_n$ consider the subtree
$T'_a\subset T'$ consisting of points fixed by a power of $a$. First
note that these subtrees are finite. Indeed, if $e$ is an edge in
$T'_a$ whose stabilizer is $<a^k>$, then exactly $k$ edges in $T'_a$
(namely, the translates of $e$ by $a$) can project to the same edge in
the quotient. By $T_a\subset T$ denote the image of $T'_a$ under the
natural map $\pi':T'\to T$ (see section 4.4). Then $T_a$ is a finite
tree (by ``local injectivity'' -- see section 4.4). Moreover, by the
equivariance of $\pi$, $T_a$ is fixed pointwise by a power of $a$.  It
then follows that if $a$ and $b$ are primitive elements with $a\neq
b^{\pm 1}$, then $T_a$ and $T_b$ can intersect in at most a point. Since
there are only finitely many orbits of the $T_a$'s, the claim follows.

It remains to rule out the possibility that
$rank(Stab(v))=\infty$. Choose a free factor $H$ of $Stab(v)$ with
$n<rank(H)<\infty$. Let $\tilde X\to T$ be a resolution such that $H$
is in the image of $\pi_1(D)\to \pi_1(X)$ for a complementary
component $D$ that corresponds to the orbit of $v$. Then
$im[\pi_1(D)\to \pi_1(X)]$ is contained in $Stab(v)$ and contains $H$,
so that its rank is $>n$. Now analyze the components of $X$ in a
similar way as above to reach a contradiction.
\end{proof}

For more details see Gaboriau-Levitt \cite{gl:rank}. They also bound the
number of orbits of branch points and their ``valences'' for small actions of
$F_n$.

\subsection{The topology of the boundary of a word-hyperbolic group}

Let $G$ be a word-hyperbolic group and $\partial G$ its boundary. The
following theorem was the motivating goal of \cite{bm:boundary}.

\begin{thm} 
If $G$ has one end, then $\partial G$ is connected and locally connected.
\end{thm}

The first part of the conclusion (that $\partial G$ is connected) was proved
in \cite{bm:boundary}, but the second was proved only under the assumption
that $\partial G$ contains no cut points. The theory of $\R$-trees was used
to establish:

\begin{thm} [Bowditch, Swarup] If $G$ has one end, then $\partial G$ contains
no cut points.
\end{thm}

\begin{proof} [Sketch of proof] For every compact metric space $M$, Bowditch
\cite{bo:dendrite} constructs a canonical map $M\to D$ to a dendrite
$D$. A compact metric space is a {\it dendrite} if it is locally
connected and each pair of points $x,y$ is joined by a unique arc,
denoted $[x,y]$. This is done as follows. Two points $x,y\in M$ are
NOT equivalent if there is a collection $C$ of cut points in $M$ that
each separate $x$ from $y$ and which is order-isomorphic to the
rationals. Bowditch argues that the quotient space $D$ is a
dendrite. Apply this construction to $M=\partial G$. Since $G$ acts on
$\partial G$, there is an induced action of $G$ on $D$. If $\partial
G$ has a cut point, then it contains a lot of cut points (translates
of the original), and Bowditch argues that $D$ is not a
point. Further, he shows that the action of $G$ on $T=D\setminus
\{endpoints\}$ has trivial arc stabilizers and is non-nesting, in the
sense that if $J$ is an arc in $T$ and $g(J)\subseteq J$, then
$g(J)=J$ (and hence $g=1$). The tree $T$ is homeomorphic to an \rt,
but there is no reason why there should be an equivariant \rt\ metric
on $T$. If there were, we could apply Theorem \ref{useful} and
conclude that $G$ splits over a 2-ended group. This is where
Sacksteder's theorem comes in. We can construct a resolution $\tilde
X\to T$ as before, but the lamination on $X$ will not have a
transverse measure. Theorem \ref{sack} provides a transverse measure
(perhaps not of full support). It is easy to see that the arc
stabilizers of the dual \rt\ are trivial. Thus $G$ splits over a
2-ended group.

The proof was completed by Swarup \cite{swarup:cut}. The idea is to
keep refining the splitting as in the proof of Theorem \ref{scott}. So
suppose inductively that $\cal G$ is a graph of groups decomposition
of $G$ with 2-ended edge groups. If $E$ is an edge group, the
endpoints of the axis of an element of $E$ are identified in the
dendrite $D$ \cite{bo:cut}. Each vertex group is word-hyperbolic and
it is quasi-convex in $G$. Let $\Lambda(V)$ denote the limit set of a
vertex group $V$ of $\cal G$. It can be argued \cite{bo:cut} that for
at least one vertex group $V$ the image of $\Lambda(V)$ in $D$ is not
a single point. It follows that the induced action of $V$ on $T$ is
nontrivial, has trivial arc stabilizers, and all edge groups contained
in $V$ are elliptic. Now apply Sacksteder's theorem again to replace
$T$ by an \rt\ $T'$ on which $V$ acts nontrivially by isometries and
with trivial arc stabilizers. The important point is that it can be
arranged that the edge groups in $V$ remain elliptic in $T'$. We then
obtain a nontrivial splitting of $V$ over two-ended groups that can be
used to refine the graph of groups decomposition $\cal G$. The final
contradiction comes from the generalized accessibility theorem
\cite{bf:bounding} that produces an upper bound to the complexity of a
(reduced) graph of groups decomposition of $G$ over two-ended groups.
\end{proof}


\providecommand{\bysame}{\leavevmode\hbox to3em{\hrulefill}\thinspace}

\end{document}